\documentclass[a4ja,10pt,leqno]{article}
\newtheorem{thm}{Theorem}[section]
\newtheorem{prop}[thm]{Proposition}
\newtheorem{cor}[thm]{Corollary}
\newtheorem{lem}[thm]{Lemma}

\newtheorem{pf}{Proof}

\newtheorem{defn}[thm]{Definition}
\usepackage{amssymb}
\usepackage{latexsym}

\newcommand{\dbar}{d\!\!\lower-0.8ex\hbox{-}}

\title{\bf Automorphisms of the Weyl manifold
\\}
\par\bigskip
\author{Naoya MIYAZAKI\thanks{The author's research is 
supported by Grant-in-Aid 
for Scientific Research (\#17540096, \#18540093), 
Ministry of Education, Culture, Sports, 
Science and Technology, Japan. 
It is also supported by Keio Gijuku Academic Funds.}
\\
Department of Mathematics, \\
Faculty of Economics, 
Keio University, \\
Yokohama, 223-8521, JAPAN.\\
\texttt{miyazaki@hc.cc.keio.ac.jp}
}

\date{}

\begin{document}
\maketitle
\par\medskip
{\abstract Assume that $M$ is a smooth manifold with a 
symplectic structure $\omega$. 
Then Weyl manifolds on the symplectic manifold $M$ are Weyl algebra bundles 
endowed with suitable transition functions. 
From the geometrical point of view, Weyl manifolds can be regarded as 
geometrizations of star products attached to $(M,\omega)$. 
In the present paper, we are concerned with 
the automorphisms of the Weyl manifold corresponding to Poincar\'{e}-Cartan 
class \footnote{$c_0$ is a $\check{\rm C}$ech 
cocycle corresponding to the symplectic structure $\omega$.} 
$[c_0+\sum_{\ell=1}^\infty c_{\ell} \nu^{2\ell}]\in \check{H}^2
(M)[[\nu^2]]$. 
We also construct modified contact Weyl diffeomorphisms\footnote
{In this paper, we often call them 
{\em lifts of symplectic diffeomorphisms}. } 
corresponding to 
symplectic diffeomorphisms of the base symplectic manifold. }
\section{Introduction}
It is well known that the concept of Lie group has a long history. 
It originated from Sophus Lie who initiated 
the systematic investigation 
of group germs of continuous transformations.  
As can be seen in introduction of a monograph 
by H. Omori \cite{o}, S. Lie seemed 
to be motivated by the followings: 
\begin{itemize}
\item To construct a theory for differential 
equation similar to Galois theory. 
\item To investigate groups 
such as continuous transformations 
that leave various geometrical structure invariant.  
\end{itemize} 
It is well known that the theory of Lie groups has expanded in two directions: 
\begin{enumerate}
\item[(A)] theory of finite-dimensional Lie groups and Lie algebras, 
\item[(B)] theory expanded to include Banach-Lie groups 
and diffeomorphisms 
those elements leave various geometrical structure invariant. 
\end{enumerate}
There are a large number of works from standpoint of (A). 
With respect to (B), there are also numerous 
works which are concerned with Banach-Lie groups and their geometrical and 
topological properties (cf. \cite{pressley-segal}). 
On the other hand, as to groups of diffeomorphisms, 
it was already known in \cite{od} that a Banach-Lie group acting effectively 
on a finite-dimensional smooth manifold is necessarily finite-dimensional. 
So there is no way to model a group of diffeomorphisms on Banach spaces 
as a manifold. Under the situation above, 
in the end of 1960s, Omori 
established 
theory of infinite-dimensional Lie groups 
called ``ILB-Lie groups"\footnote
{``ILB" means {\em inductive limit of Banach spaces}. }  
beyond Banach-Lie groups, taking ILB-chains as model spaces
in order to treat groups of diffeomorpshisms on a manifold 
(see \cite{o} for the precise definition). 
Shortly after his works, Omori et al. \cite{omyk} 
introduced the 
definition of Lie group modeled on a Fr\'{e}chet space equipped with 
a certain property called  ``regurality"  by relaxing the conditions of 
ILB-Lie group. Roughly speaking, regularity means that the smooth curves in 
the Lie algebra integrate to smooth curves in the Lie group in a smooth way
(see also \cite{mi}, \cite{o} and \cite{omyk1}). 
Using this notion, they studied subgroups of a group of 
diffeomorphisms, and 
the group of invertible Fourier integral operators with 
suitable amplitude functions on a manifold. 
For technical reasons, they assumed that the base manifold is compact 
(cf. \cite{ars}, \cite{ars1}, \cite{ars2}, \cite{mic} and \cite{omyk}).  Beyond a compact base manifold, 
in order to treat groups of diffeomorphisms on 
a {noncompact} manifold, we need more general category of Lie groups, i.e. 
infinite-dimensional Lie groups modeled on 
locally convex spaces which are Mackey complete 
(see \S 2. See also \cite{er} and \cite{km}). 

In this paper, 
we are concerned with the group
\footnote{See Definition~\ref{MCWD} for the precise definition.} 
${\rm Aut}(M,*)$ of all modified contact Weyl diffeomorphisms 
on a contact Weyl manifold 
over a symplectic manifold $(M,\omega)$, 
where a contact Weyl manifold introduced by A. Yoshioka in \cite{y} 
is a geometric realization of star product introduced in \cite{bffls}. 
In this context, a modified contact Weyl diffeomorphism 
is regarded as an automorphism 
of star product. 
As to the group 
${\rm Aut}(M,*)$, we have the following. 
\begin{thm}\label{main-theorem-omori-proceedings}
\begin{enumerate}
\item 
Set $$\underline{\rm Aut}(M,*)=\{\Phi\in {\rm Aut}(M,*)\,|\, 
\Phi\mbox{ induces the base identity map.}\}.$$  
Then $\underline{\rm Aut}(M,*)$  is 
a Lie group modeled on a Mackey complete locally convex space. 
\item Any element $\Psi\in {\rm Aut}(M,*)$ induces a 
symplectic diffeomorphism on the base manifold and there exists 
a group homomorphism $p$ from 
${\rm Aut}(M,*)$ onto ${\rm Diff}(M,\omega)$, where ${\rm Diff}(M,\omega)$ 
is the regular Lie group of 
all symplectic diffeomorphisms on the symplectic manifold $(M,\omega)$. 
\item 
The group ${\rm Aut}(M,*)$ is a Lie group modeled 
on a Mackey complete locally convex space
\footnote{If the base manifold is compact, the model spaces of 
${\rm Aut}(M,*)$ and $~\underline{\rm Aut}(M,*)$ are Fr\'{e}chet spaces.}. 
\item Under the same assumption above,  
$$1\rightarrow \underline{\rm Aut}(M,*)\rightarrow {\rm Aut}(M,*)  
\rightarrow {\rm Diff}(M,\omega)\rightarrow 1$$ 
is a short exact sequence of Lie groups. 
\item The groups 
$\underline{\rm Aut}(M,*)$ and ${\rm Aut}(M,*)$ are regular Lie groups. 
\end{enumerate}
\end{thm}
\noindent
We note that 
this result in formal deformation quantization 
might be regarded as a counterpart 
of the result of the regular Lie group sturcture for 
the group $G{\mathcal F}_0(N)$ of invertible Fourier integral operators 
with classical symbols of 
order 0 on a riemannian manifold $N$ in micro-local analysis
(see \cite{ars2}, \cite{er} and \cite{omyk}). 
Moreover it is also known that the following sequence 
\begin{equation}
1\rightarrow G\Psi_0(N) 
\rightarrow G{\mathcal F}_0(N) \rightarrow {\rm Diff}(S^*N,\theta) 
\rightarrow 1  
\end{equation}
is exact, where $G\Psi_0(N)$ (resp. ${\rm Diff}(S^*N,\theta)$) denotes 
the group of invertible pseudo-differential operators with 
classical symbols of order 0 (resp. the group of contact diffeomorphisms on 
the unit cosphere bundle $S^*N$ with the contact structure $\theta$). 
\par\medskip
Remark that from the point of view of differential geometry, 
a contact Weyl manifold might be seen as a ``{\it prequantum bundle}" 
over a symplectic manifold $(M,\omega)$ 
where the symplectic structure $\omega$ is not necessarily {\it integral}, 
and a modified contact Weyl diffeomorphim can be regarded 
as a quantum symplectic diffeomorphism over a ``prequantum bundle". 

As is well known, theory of infinite-dimensional 
Lie algebras including Kac-Moody algebras has made rapid and remarkable 
progress for the past two decades involving completely integrable system 
(Sato's theory), loop groups, conformal field theory and quantum groups. 
However, it would be rather difficult for me to review this fruitful field. 
A definite treatment of the infinite-dimensional Lie algebras is 
found in Kac~\cite{kac}, Tanisaki~\cite{tanisaki} and 
Wakimoto~\cite{wakimoto}. 
\par\medskip

\section{Infinite-dimensional Lie groups}
In this section we give a survey of regular Lie groups. 
For the purpose, we first recall Mackey completeness, 
see the excellent monographs \cite{j}, \cite{km} for details.  
\begin{defn}
A locally convex space $E$ 
is called a {\rm Mackey complete} 
$(${\rm MC} for short$)$ if one of the following equivalent 
conditions is satisfied: 
\begin{enumerate}\label{completeness}
\item For any smooth curve $c$ in $E$ there is a smooth curve $C$ in $E$ 
with $C'=c$. 
\item If $c:{\mathbf R}\rightarrow E $ is a curve 
such that $l\circ c:{\mathbf R}\rightarrow {\mathbf R}$ is smooth for all 
$\ell \in E^*$, then $c$ is smooth. 
\item Locally completeness: For 
every absolutely convex closed bounded\footnote{
A subset $B$ is called bounded if it is absorbed 
by every $0$-neighborhood in $E$, 
i.e. for every $0$-neighborhood ${\mathcal U}$, there exists 
a positive number $p$ 
such that $[0,p]\cdot B\subset {\mathcal U}$.} subset $B$, 
$E_B$ is complete, where $E_B$ is 
a normed space linearly generated by $B$ with 
a norm $p_B(v)=\inf\{\lambda>0|v\in \lambda B\}$.
\item Mackey completeness: any Mackey-Cauchy net converges in $E$.
\item Sequential Mackey completeness: any Mackey-Cauchy  
sequence converges in $E$.  
\end{enumerate}
Here 
a net $\{x_\gamma\}_{\gamma\in \Gamma}$ is called {\rm Mackey-Cauchy} 
if there exists a bounded set $B$ and  
a net $\{\mu_{\gamma,\gamma'}\}_{(\gamma,\gamma')\in \Gamma\times \Gamma}$ 
in $\mathbf R$ converging to $0$, 
such that 
$x_\gamma-x_{\gamma'} \in  \mu_{\gamma,\gamma'}B = 
\{\mu_{\gamma,\gamma'}\cdot x|x\in B \}$.
\end{defn}
\par\medskip
We recall 
the fundamentals relating to infinite-dimensional differential geometry. 
\begin{enumerate}
\item Infinite-dimensional manifolds 
are defined on Mackey complete locally convex spaces in much the same way 
as ordinary manifolds are defined on finite-dimensional 
spaces. 
In this paper, 
a manifold equipped with a smooth group operation is referred to 
as a Lie group. 
Remark that in the category of infinite-dimensional Lie groups, 
the existence of exponential maps is not ensured 
in general, and even if an exponential map exists, 
the local surjectivity of it does not hold  
(cf. Definition \ref{regular}). 
\item A {\it kinematic tangent vector} (a tangent vector for short) 
with a foot point $x$ of an infinite-dimensional manifold $X$ modeled on 
a Mackey complete locally convex space $F$ is a pair $(x,X)$ with 
$X\in F$, and let $T_xF=F$ be the space of all tangent vectors with foot 
point $x$. It consists of all derivatives $c'(0)$ at $0$ of smooth curve 
$c:{\mathbf R} \rightarrow F$ with $c(0)=x$. 
Remark that operational tangent vectors viewed as derivations  
and kinematic tangent vectors via curves differ in general.  
A kinematic vector field is a smooth section of kinematic vector bundle 
$TM\rightarrow M$. 
\item We set $\Omega^k(M)=C^\infty(L_{skew}
({TM\times\cdots\times TM},M\times{\mathbf R}))$ and call it the space of 
{\it kinematic differential forms}, where ``skew" denotes ``skew-symmetric". 
Remark that the space of kinematic  
differential forms turns out to be the right ones for 
calculus on manifolds; especially for 
them the theorem of de~Rham is proved. 
\end{enumerate}
Next we give the precise definition of regularity (cf. \cite{mi}, 
\cite{o}, \cite{omyk} and \cite{omyk1}): 
\begin{defn}\label{regular}
A Lie group $G$ modeled on a Mackey complete locally convex space 
$\frak{G}$ is called a {\rm regular} 
Lie group if one of the following equivalent 
conditions is satisfied 
\begin{enumerate}
\item For each $X \in C^\infty({\mathbf R}, \frak{G})$, there exists 
$g\in C^\infty({\mathbf R},G)$ satisfying 
\begin{eqnarray}
&&g(0)=e, \quad \frac{\partial}{\partial t}g(t)=R_{g(t)}(X(t)),
\label{regularity-eq}
\end{eqnarray} 
\item For each $X \in C^\infty({\mathbf R}, \frak{G})$, there exists 
$g\in C^\infty({\mathbf R},G)$ satisfying 
\begin{eqnarray}
&&g(0)=e, \quad \frac{\partial}{\partial t}g(t)=L_{g(t)}(X(t)),
\end{eqnarray} 
\end{enumerate} 
where $R(X)$ {\rm (}resp. $L(X)${\rm )} is the right 
{\rm (}resp. left{\rm )} invariant vector field   
defined by the right{\rm (}resp. left{\rm )}-translation of a 
tangent vector $X$ at $e$. 
\end{defn}
The following lemma is useful (cf. \cite{km}, \cite{mi}, \cite{omyk} and 
\cite{omyk1}): 
\begin{lem}\label{short-exact-sequence}
Assume that 
\begin{equation}
1 \rightarrow N \rightarrow G \rightarrow H \rightarrow 1
\end{equation}
is a short exact sequence of Lie groups with a 
local smooth section\footnote{Remark that this does not give 
global splitting of the short exact sequence.} $j$ from a neighborhood 
$U \subset H$ of $1_H$ into $G$, and 
$N$ and $H$ are regular. 
Then 
$G$ is also regular. 
\end{lem}
To end this section, we remark that the fundamental properties of 
principal regular Lie group bundle $(P,G)$ over $M$. 
Note that these properties are ordinary properties 
for principal finite-dimensional Lie group bundles.    
\begin{enumerate}
\item The parallel transformation is well defined. 
\item The horizontal distribution ${\mathcal H}$ of a flat connection is 
integrable, i.e. 
there exists an integral submanifold for ${\mathcal H}$ at each point. 
\end{enumerate}
\section{Deformation Quantization}
Mathematically the concept of quantization 
originated from H. Weyl \cite{weyl}, who introduced a map 
from classical observables 
(functions on the phase space) 
to quantum obsevables (operators on  Hilbert space). 
The inverse map 
was constructed by E. Wigner 
by interpreting functions 
(classical observables) as symbols of operators. 
It is known that 
the exponent of the bidifferential operator (Poisson bivector) coincides with 
the product formula of Weyl type symbol calculus 
developed by L. H\"{o}rmander who established the theory 
of pseudo-differential operators and used them to study 
partial differential equations (cf. \cite{ku} and \cite{moyal}). 

In the 1970s, 
supported by the mathematical developments above, 
Bayen, Flato, Fronsdal, Lichnerowicz and Sternheimer \cite{bffls} 
considered quantization as a deformation of 
the usual commutative product of classical observables 
into a noncommutative associative product which is 
parametrized by the Planck constant $\hbar$ and 
satisfies the correspondence principle.  
Nowadays 
deformation quantization, or moreprecisely, star product
has gained support from geometricians and mathematical physicists. 
In fact, it plays an important role to give passage from 
Poisson algebras of classical observables to 
noncommutative associative algebras of quantum observables. 
In the approach above, the precise definition of 
the space of quantum observables and star product 
is given in the following way(cf. \cite{bffls}):   
\begin{defn}{}\label{FDQ}
A 
{\rm star product} of
Poisson manifold $(M,\pi)$
is a product $*$
on the space $C^{\infty}(M)[[\hbar]]$ 
of formal power series of parameter
$\hbar$ with coefficients in $C^{\infty}(M)$,
defined by
\begin{eqnarray}
f * g \!\!&=&\!\!
 fg+\hbar\pi_1 (f,g)+\cdots\nonumber
+\hbar^n \pi_n (f,g)+\cdots, ~~\forall f,g\in C^{\infty}(M)[[\hbar]]
\nonumber
\end{eqnarray}
satisfying \\
{\rm (a)} $*$ is associative, \\
{\rm (b)} $\pi_1 (f,g)=\frac{1}{2\sqrt{-1}}\{f,g\},$ \\
{\rm (c)} each $\pi_n~(n\ge 1)$ is a ${\mathbf C}[[\hbar]]$-bilinear and
bidifferential operator, \\
where $\{,\}$
is the Poisson bracket defined by the Poisson structure $\pi$.
\end{defn}
A deformed algebra (resp. a deformed algebra structure) 
is called a {\it star algebra} (resp. a {\it star product}).
Note that on a symplectic vector space  ${\mathbf R}^{2n}$,
there exists the ``canonical'' deformation 
quantization, the so-called {Moyal product}: 
$$f*g=f\exp\bigr[ \frac{\nu}{2}\stackrel{\leftarrow}{\partial_x}\land 
                 \stackrel{\rightarrow}{\partial_y} \bigl] g, $$
where $f,~g$ are smooth functions of a Darboux coordinate 
$(x,y)$ on ${\mathbf R}^{2n}$ 
and $\nu=i\hbar$. 

The existence and classification problems of star products
have been solved 
by succesive steps from special classes of symplectic
manifolds to general Poisson manifolds. 
Because of its physical origin and motivation, the problems 
of deformation quntization was first considered for symplectic manifolds, 
however, the problem of deformation quantization is naturally 
formulated for the Poisson manifolds as well. 
For example, Etingof and Kazhdan proved every 
Poisson-Lie group can be quantized in the sense above, and investigated  
quantum groups as deformation quantization of Poisson-Lie groups. 
After their works, for a while,  
there were no specific developments for existence problems 
of deformation quantization on any Poisson manifold. 
The situation drastically changed when M. Kontsevich \cite{k} proved 
his celebrated formality theorem. As a collorary, he showed that 
deformation quantization exists on 
any Poisson manifold. 
(cf. \cite{dl}, \cite{f}, \cite{k}, \cite{ommy}, \cite{omy}, \cite{s} and 
\cite{y}).  
\section{Weyl manifold and contact Weyl manifold over a symplectic manifold}
\subsection{Definition of (contact) Weyl manifold}
As mentioned in the introduction, 
by Omori-Maeda-Yoshioka, 
for a symplectic manifold, the notion of Weyl manifold was introduced. 
Later, Yoshioka \cite{y} proposed the notion of contact Weyl manifold as 
a bridge joining the theory of Weyl 
manifold and the Fedosov approach to quantization. 
In order to recall the construction of a contact Weyl manifold, 
we have to give the precise definitions of fundamental algebras. 
\begin{defn}\label{Weyl-alg}
\begin{enumerate}
\item An associative algebra $W$ is called a {\rm Weyl algebra} if $W$ is 
formally generated by  
$\nu,Z^1,\ldots,Z^n,Z^{n+1},\ldots,Z^{2n}$
satisfying the following commutation relations:
\begin{equation}\label{CCR}
[Z^i,Z^j]=\nu\Lambda^{ij},~~[\nu,Z^i]=0,
\end{equation}
where $\scriptsize{\Lambda=\left[\begin{array}{cc} 0 & -1_n\\ 1_n &0
\end{array}\right]}$, and the product of this algebra is denoted by $*$.  
This algebra has the canonical involution 
$\bar{}$ such that 
\begin{equation}\label{involution}
\overline{a*b}=\bar{b}*\bar{a},\quad \bar{\nu}=-\nu, ~~~\bar{Z^i}=Z^i. 
\end{equation}
We also define the degree $d$ by $d(\nu^lZ^\alpha)=2l+|\alpha|.$ 
\item A Lie algebra $C$ is called a {\rm contact Weyl algebra}  
if $C=\tau {\mathbf C}\oplus W$ 
with an additional generator $\tau$ 
satisfying the following relations: 
\begin{equation}\label{CW-relation}
[\tau,\nu]=2\nu^2,~~~[\tau , Z^i]=\nu Z^i,
\end{equation} 
and $\bar{}$ is naturally extended by $\bar{\tau}=\tau$.
\end{enumerate} 
\end{defn}
Remark that the relation (\ref{CCR}) is nothing but 
the commutation relation of the Moyal product, and called 
the {\it canonical commutation relation}. 
It is well known that the ordering problem appears 
when we realize this algebra explicitly. 
In this paper, we mainly 
use the Weyl ordering (the Moyal product). 
See Appendix~6.2, for related topics. 
\par\medskip\noindent
\begin{defn}\label{nu-auto}
\begin{enumerate}
\item A ${\mathbf C}[[\nu]]$-linear isomorphism $\Phi$ from $W$ onto $W$ 
is called a {\rm $\nu$-automorphism} of Weyl algebra $W$ if 
   \begin{enumerate}
   \item[]{\rm (a)} $\Phi(\nu)=\nu $,
   \item[]{\rm (b)} $\Phi(a*b)=\Phi(a)*\Phi(b)$, 
   \item[]{\rm (c)} $\Phi(\bar{a})=\overline{\Phi(a)}.$
   \end{enumerate}
\item A ${\mathbf C}[[\nu]]$-linear isomorphism $\Psi$ from $C$ onto $C$ 
is called a {\rm $\nu$-automorphism} of contact Weyl algebra $C$ if 
   \begin{enumerate}
   \item[]{\rm (a)} $\Psi$ is an algebra isomorphism, 
   \item[]{\rm (b)} $\Psi|_W$ is a $\nu$-automorphism of Weyl algebra.
   \end{enumerate}
\end{enumerate}
\end{defn}
In order to explain the construction of contact Weyl manifolds, it is 
useful to recall how to construct 
prequantum line bundles, which play an crucial role in 
the theory of Souriau-Kostant (geometric) quantization\cite{w}. 
This bundle is constructed in the following way: 
Let $\omega $ be an integral symplectic structure, 
then we have  
$d(\theta_\alpha)=(\delta\omega)_{\alpha},~
d(f_{\alpha\beta})= 
(\delta\theta)_{\alpha\beta},~
c_{\alpha\beta\gamma}=
(\delta f)_{\alpha\beta\gamma}
$ 
where ${\mathcal U}=\{U_\alpha\}$ is a good covering of a symplectic manifold 
$(M,\omega)$, 
$f_{\alpha\beta}$ 
(resp. $\theta_\alpha$) is a local function (resp. a local 1-form) defined on 
an open set 
$U_\alpha\cap U_\beta$ (resp.~$U_\alpha$),   
$d$ is the deRham exterior differential operator, and 
$\delta$ is the $\check{\rm C}$ech coboundary operator. 
Setting $h_{\alpha\beta}=\exp [2\pi i f_{\alpha\beta}]$, 
we see 
that 
\begin{equation}\label{bundle-condition}
\theta_\alpha-\theta_\beta=\frac{1}{2\pi i}d \log h_{\alpha\beta}.  
\end{equation}
This equation ensures the exsistence of a line bundle defined by   
\begin{eqnarray}
&&L=
\coprod (U_\alpha\times {\mathbf C})/\stackrel{h_{\alpha\beta}}{\thicksim} , 
\qquad\nabla_\xi(\phi_\alpha 1_\alpha)
     =(\xi\phi_\alpha + 2\pi i \theta_\alpha(\xi)\phi_\alpha)1_\alpha .   
\end{eqnarray}
This gives the desired bundle with a connection whose curvature equals 
$\omega$. 
\par
Inspired by the idea above, Yoshioka 
proposed the notion of contact Weyl manifold and 
obtained the fundamental results (cf. \cite{y}). 
To state the precise definition of contact Weyl manifold 
and theorems related to them, we recall the definitions 
of {\em Weyl continuation} and 
{\em locally modified contact Weyl diffeomorphism}:  
\begin{defn}
Set $(X^1,\ldots,X^n,Y^1,\ldots, Y^n):=(Z^1,\ldots,Z^n,Z^{n+1},\ldots,
Z^{2n})$ {\rm (}see Definition~\ref{Weyl-alg}{\rm )}.   
Consider the trivial contact Weyl algebra bundle $C_U:=U\times C$ over a 
local Darboux chart $(U;({x},{y}))$. 
A section 
\begin{eqnarray}
\label{weyl-conti}
\nonumber
f^{\#}:= f({x}+X,{y}+Y)
= \sum_{\alpha\beta}
\frac{1}{\alpha!\beta!}\partial_x^\alpha\partial_y^\beta  
f({x},{y})X^\alpha Y^\beta \in \Gamma(C_U)
\end{eqnarray} 
determined by a local smooth function $f\in C^\infty(U)$ 
is called a {\rm Weyl function}, and $\#:f\mapsto f^{\#}$ 
is referred to as {\rm Weyl continuation}. 
We denote by ${\mathcal F}_U$ the set of all Weyl functions on $U$.  

A bundle map $\Phi:C_U\rightarrow C_U$ is 
referred to as a  {\rm locally modified contact Weyl diffeomorphism} if 
it is a fiberwise $\nu$-automorphism of the contact Weyl algebra 
and its pull-back preserves the set of all Weyl functions ${\mathcal F}_U.$ 
\end{defn}

\begin{defn}
Let $\pi:C_M\rightarrow M$ be a 
locally trivial bundle with a fiber being isomorphic to the 
contact Weyl algebra over a symplectic manifold $M$.
Take an atlas $\{(V_\alpha,\varphi_\alpha)\}_{\alpha\in A}$ of 
$M$ such that
$\varphi_\alpha:V_\alpha\rightarrow U_\alpha\subset {\mathbf R}^{2n}$ 
gives a local Darboux coordinate for every $\alpha\in A$. 
Denote by $\Psi_\alpha:C_{V_\alpha}\rightarrow C_{U_\alpha}$ a 
local trivialization and by 
$\Psi_{\alpha\beta}=\Psi_\beta\Psi_\alpha^{-1}:
C_{U_{\alpha\beta}}\rightarrow C_{U_{\beta\alpha}}$ the glueing map, 
where 
$C_{V_\alpha}:=\pi^{-1}(V_\alpha)$, 
$U_{\alpha\beta}:=\varphi_\alpha(V_\alpha\cap V_\beta), 
~U_{\beta\alpha}:=\varphi_\beta(V_\alpha\cap V_\beta),~ 
C_{U_{\alpha\beta}}:=\Psi_\alpha(C_{V_\alpha}|_{V_\alpha\cap V_\beta})$, etc. 
Under the notations above, 
\begin{equation}
\label{cwm}
\bigl(\pi: C_M\rightarrow M,  
\{\Psi_\alpha:C_{V_\alpha}\rightarrow C_{U_\alpha}\}_{\alpha\in A}\bigl)
\end{equation}  
is called a {\rm contact Weyl manifold}\footnote{
Let $W_U$ be a trivial Weyl algebra bundle attached to 
a Darboux coordinate neighborhood. A bundle map $\Phi:W_U\rightarrow W_U$ is 
referred to as a  {\rm local Weyl diffeomorphism} if 
it is a fiberwise $\nu$-automorphism of the Weyl algebra 
and its pull-back preserves the set of all Weyl functions ${\mathcal F}_U.$ 
Originally, 
using the notion of local Weyl diffeomorphisms, 
Omori-Maeda-Yoshioka gave the definition of Weyl manifold. 
\begin{defn}
Let 
$\pi: W_M\rightarrow M$ be a 
locally trivial bundle with a fiber being isomorphic to the 
Weyl algebra over a symplectic manifold $M$.
Take an atlas $\{(V_\alpha,\varphi_\alpha)\}_{\alpha\in A}$ of 
$M$ such that
$\varphi_\alpha:V_\alpha\rightarrow U_\alpha\subset {\mathbf R}^{2n}$ 
gives a local Darboux coordinate for every $\alpha\in A$. 
Denote by $\Phi_\alpha:W_{V_\alpha}\rightarrow W_{U_\alpha}$ a 
local trivialization and by 
$\Phi_{\alpha\beta}=\Phi_\beta\Phi_\alpha^{-1}:
W_{U_{\alpha\beta}}\rightarrow W_{U_{\beta\alpha}}$ the glueing map, 
where 
$W_{V_\alpha}:=\pi^{-1}(V_\alpha)$, 
$U_{\alpha\beta}:=\varphi_\alpha(V_\alpha\cap V_\beta), 
~U_{\beta\alpha}:=\varphi_\beta(V_\alpha\cap V_\beta),~ 
W_{U_{\alpha\beta}}:=\Phi_\alpha(W_{V_\alpha}|_{V_\alpha\cap V_\beta})$, etc. 
Under the notations above, 
\begin{equation}
\label{wm}
\bigl(\pi: W_M\rightarrow M,  
\{\Phi_\alpha:W_{V_\alpha}\rightarrow W_{U_\alpha}\}_{\alpha\in A}\bigl)
\end{equation}  
is called a {\rm Weyl manifold} 
if the glueing maps $\Phi_{\alpha\beta}$ 
are local Weyl diffeomorphisms. 
\end{defn}
} 
if the glueing maps $\Psi_{\alpha\beta}$ 
are modified contact Weyl diffeomorphisms. 
\label{cwm}
\end{defn}
\begin{thm}\label{y}
Let $(M,\omega)$ be an arbitrary {\rm (}not necessarily integral{\rm )} 
symplectic manifold.  
For any closed form 
$\Omega_M(\nu^2)=\omega+\omega_2\nu^2+\omega_4\nu^4+\cdots$, 
where $\nu=\sqrt{-1}\hbar$ is a formal parameter, 
there exists a contact Weyl manifold 
$C_M$ with a connection 
$\nabla^Q$ whose curvature 
equals $ad[\frac{1}{\nu}\Omega_M(\nu^2)]$, and the restriction of 
$\nabla^Q$ to $W_M$ is flat, where
$W_M$ is the Weyl algebra bundle 
associated to $M$ equipped with the canonical fiber-wise 
product $\hat{*}$. 
\end{thm}
This bundle $C_M$ is called a contact Weyl manifold 
equipped with a {\it quantum connection} $\nabla^Q$. 
Yoshioka \cite{y} also proved that 
the connection $\nabla^Q|_{W_M}$ is essentially the same 
as a Fedosov connection 
$\nabla^{W}$ \cite{f}. 
It is known (cf. \cite{y} and \cite{ommy}) that 
\begin{thm}\label{theorem-bijection}
There is a bijection between 
the space of the isomorphism classes of contact Weyl manifolds 
with quantum connections
and 
$[\omega]+\nu^2 H_{dR}^2(M)[[\nu^2]]$, 
which assigns a class $[\Omega_M(\nu^2)]=[\omega+\omega_2\nu^2+\cdots]$ to a 
contact Weyl manifold with quantum connection 
$(C_M\rightarrow M,\{\Psi_\alpha\},\nabla^Q)$. 
\end{thm}
\begin{pf}
It is already known that there is a bijection between 
the space of the isomorphism classes of Weyl manifolds 
and $[\omega]+\nu^2 H_{dR}^2(M)[[\nu^2]]$ 
$($cf. subsection~\ref{subsection-4-2}$)$. 
Generalizing straightly the proof of this result, we can prove 
Theorem~\ref{theorem-bijection}. 
$\square$
\end{pf}
The flatness of $\nabla^Q|_{W_M}$ ensures the existence of 
a linear isomorphism $\#$ 
between 
$C^\infty(M)[[\nu]]$ and 
${\mathcal F}_M$ the space of all parallel sections with respect to 
the quantum connection restricted to $W_M$. 
An element of ${\mathcal F}_M$ is called a {\it Weyl function}.  
Using this map $\#$, we can recapture a star product in the following way: 
\begin{equation}\label{recap}
f*g=\#^{-1}\bigl(\#  (f) \hat{*} \# (g)\bigl). 
\end{equation}    
Furthermore, it is known that the following (cf. \cite{ommy}, see also 
\cite{de} and \cite{gr}):
\begin{thm}\label{class}
There is a bijection between 
the space of the equivalence classes of star products 
and
$[\omega]+\nu^2 H_{dR}^2(M)[[\nu^2]]$.
\end{thm}

\subsection{Poincar\'{e}-Cartan classes (Deligne relative classes)}
\label{subsection-4-2} 
We begin this subsection with the fundamental facts and definitions. 
Set $\tilde{\tau}_U=\tau+\sum z^i\omega_{ij}Z^j$ 
where $U\subset {\mathbf R}^{2n}$ is an open subset and 
$\omega_{ij}dz^i\land dz^j$ stands for the symplectic structure. 
Then for any modified contact Weyl diffeomorphism, we may set ${\Psi}|_{{ C}_U}^\ast(\tilde{\tau}_U)=
a\tilde{\tau}_U+F$, where $a\in C^\infty(U),\, 
F\in\Gamma({ W}_U)$,  
where $W_U$ is a trivial bundle 
$W_U=U\times W$. 
Furthermore it is known that the following (cf. Lemma 2.21 in \cite{y}).
\begin{prop}\label{representative1} 
Let $U$ be an open set in ${\mathbf R}^{2n}$, $\Psi$ 
a modified contact Weyl diffeomorphism 
and $\phi$ the induced map on the base manifold.
Then the pullback of $\tilde{\tau}_{\phi(U)}$ by $\Psi$ can be written as 
\begin{equation}
\Psi^\ast\tilde{\tau}_{\phi(U)}=\tilde{\tau}_U+f^{\#}+a(\nu^2),
\end{equation}
for some Weyl functions $f^{\#}:=\#(f)\in {\mathcal F}_U$
 with $\bar{f^\#}=f^\#$ and 
$a(\nu^2)\in C^\infty(U)[[\nu^2]]$.   
\end{prop}
\begin{defn} 
A modified contact Weyl diffeomorphism ${\Psi}$ 
is called a {\rm contact Weyl diffeomorphism\,(CWD, for short)} if 
\begin{equation}
{\Psi}^\ast\tilde{\tau}_{U'}=\tilde{\tau}_U+f^{\#}. 
\end{equation}
\label{defn4-10}
\end{defn}
For a contact Weyl diffeomorphism, 
we obtain the following  
(see Corollary 2.5 in \cite{ommy} and Proposition 2.24 in \cite{y}). 
\begin{prop}\label{2.4ommyrevised}
\label{ommy-a}Assume that a map $\Psi$ is a contact Weyl diffeomorphism. 
\begin{enumerate}
\item If the diffeomorphism $\varphi$ on the base map induced by 
$\Psi$ is the identity, there exists uniquely a Weyl function $g^{\#}(\nu^2)$ 
such that $$\Psi=\exp[ad(\frac{1}{\nu}g^{\#}(\nu^2))].$$
\item $\Psi|_{{W}_U}=1$ if and 
only if there exists an element $ c(\nu^2)\in {\mathbf R}[[\nu^2]]$ such that 
$$
\nonumber \Psi = \exp [\frac{1}{\nu}ad (c(\nu^2))].
$$
\end{enumerate}
\end{prop}
From this proposition, we can define the Poincar\'{e}-Cartan class 
in the following way (\cite{ommy}). 
Assume that $W_M=\{ (W_{U_\alpha},\Phi_{\alpha\beta}) \}$ is a Weyl manifold. 
Then $\Phi_{\alpha\beta}\Phi_{\beta\gamma}\Phi_{\gamma\alpha}$ is 
the identity on each $W_{U_{\alpha\beta\gamma}}$. 
According to 2 of Proposition \ref{ommy-a},  we have 
$$\Phi_{\alpha\beta}\Phi_{\beta\gamma}\Phi_{\gamma\alpha}=
\exp[\frac{1}{\nu}(c_{\alpha\beta\gamma}(\nu^2))],~
~(\exists c_{\alpha\beta\gamma}(\nu^2)\in{\mathbf R}[[\nu^2]]).$$
We can show that 
$\{c_{\alpha\beta\gamma}\}$ is a $\check{\rm C}$ech 2-cocycle, and then 
it defines a $\check{\rm C}$ech 2-class.  
\begin{defn}\label{pcdef}
We refer to this cocycle $($resp. class$)$ 
as the {\rm Poincar\'{e}-Cartan cocycle (resp. class)} 
and denote it by $\{c_{\alpha\beta\gamma}\}$ $($resp. $c(W_M)$$)$.  
\end{defn}
For the Poincar\'{e}-Cartan class we have the following. 
\begin{thm}
For any $ c=c^{(0)}+\sum_{i=1}^\infty c^{(2i)}\nu^{2i}
~(c^{(2i)}\in \check{H}^2(M;{\mathbf R}))$ 
such that $[c^{(0)}]$ corresponds to the class of symplectic 2-form, 
there exists a family of 
contact Weyl  diffeomorphisms
$ \{\Psi_{\alpha\beta}:{C}_{U_{\alpha\beta}}\rightarrow
{C}_{U_{\beta\alpha}}\}$, such that 
$\Psi_{\alpha\beta\gamma}|_{W_{U_{\alpha\beta\gamma}}}=1$, where 
$\Psi_{\alpha\beta\gamma}:=\Psi_{\alpha\beta}
\Psi_{\beta\gamma}\Psi_{\gamma\alpha}$, and  
$\{c_{\alpha\beta\gamma}(\nu^2)\}$ defines 
a $\check{C}$ech 2 cohomology class which coincides with 
$c$. Moreover there is one to one correspondence between 
the set $\frak{PC}(M)$ of Poincar\'{e}-Cartan classes 
and the set $\frak{W}(M)$ of isomorphism classes of Weyl manifolds.  
\label{PCclass} 
\end{thm}

\begin{pf}
The proof is already known in \cite{ommy}, but we here give an outline of it 
for readers.
First we show that for any cocycle $\{c_{\alpha\beta\gamma}(\nu^2)\}$, 
there exists a Weyl manifold such that 
$c({W}_M)=[c_{\alpha\beta\gamma}(\nu^2)]$. 
Suppse that $c=\sum_{k\ge 0}\nu^{2k}c^{(2k)}\in \check{H}^2(M) [[\nu^2]]$ 
is given. 
According to the existence theorem of Weyl manifold in \cite{omy}, we may 
start with a Weyl manifold ${W}_M^{(0)}$ with a 
Poincar\'{e}-Cartan cocycle $\{c_{\alpha\beta\gamma}^{(0)}\}$, and changing 
patching Weyl diffeomorphisms we construct a Weyl manifold with a 
Poincar\'{e}-Cartan class $c$. 
Let $\Phi^*_{\alpha\beta}:{\cal F}(W_{U_{\alpha\beta}})\rightarrow
{\cal F}(W_{U_{\beta\alpha}})$ be the glueing 
Weyl diffeomorphism of ${W}_M^{(0)}$ and let  
$\Psi^*_{\alpha\beta}$ be its extension as a contact Weyl diffeomorphism. 
Let $\{c_{\alpha\beta\gamma}^{(2k)}\}$ be a $\check{\rm C}$ech  
cocycle belonging to $c^{(2k)}$. 
Since the sheaf cohomology $H^2(M;{\cal E})$ 
of the sheaf of germs $C^\infty$-functions ${\cal E}$, 
there is $h_{\alpha\beta}^{(2)}\in C^\infty(U_{\alpha\beta})$ on 
each $U_{\alpha\beta}$ such that 
\begin{equation}\label{id3.6}
-c_{\alpha\beta\gamma}^{(2)}=h_{\alpha\beta}^{(2)}
+\varphi^*h_{\beta\gamma}^{(2)}+\varphi_{\alpha\gamma}^*
h_{\gamma\alpha}^{(2)} .\end{equation}
Replace $\Psi_{\alpha\beta}^*$ by 
$\dot{\Psi}_{\alpha\beta}^*
=\Psi_{\alpha\beta}^*e^{ad (\nu\tilde{h}_{\beta\alpha}) }
$ 
as glueing diffeomorphism for each 
$V_\alpha\cap V_\beta\not=\phi $. 
Then according to the formula 
$ 
\Psi_{\alpha\beta}^*  e^{ad(h)}=e^{ad(\Psi_{\alpha\beta}^*h)}
\Psi_{\alpha\beta}^*
$ 
for $h\in {\cal F}({W}_{U_{\alpha\beta}})$, 
we see 
$$\dot{\Phi}_{\alpha\beta}^*\dot{\Phi}_{\beta\gamma}^*
\dot{\Psi}_{\gamma\alpha}^*
=\Psi_{\alpha\beta}^*\Psi_{\beta\gamma}^*\Psi_{\gamma\alpha}^*
e^{ad(\nu \Psi_{\alpha\beta}^*\tilde{h}_{\beta\alpha})}
e^{ad(\nu \Psi_{\alpha\gamma}^*\tilde{h}_{\gamma\beta})}
e^{ad(\nu \Psi_{\alpha\alpha}^*\tilde{h}_{\alpha\gamma})},  
$$ 
where we set 
$\tilde{h}_{\alpha\beta}=(h_{\alpha\beta}^{(2)})^\#+\nu^2r_{\alpha\beta}^\#$ 
for a function 
$r_{\alpha\beta} \in C^\infty(U_{\alpha\beta})[[\nu^2]]$. 
 By (\ref{id3.6}),
we have 
\begin{equation}\label{id3.8}
e^{ad(\nu \Psi_{\alpha\beta}^*\tilde{h}_{\beta\alpha})}
e^{ad(\nu \Psi_{\alpha\gamma}^*\tilde{h}_{\gamma\beta})}
e^{ad(\nu \Psi_{\alpha\alpha}^*\tilde{h}_{\alpha\gamma})}
=e^{\nu^2c_{\alpha\beta\gamma}^{(2)}ad  ( \nu^{-1} ) }~~mod~\nu^4.
\end{equation}
By working on the term $\nu^4,~\nu^6,\cdots $, 
we can tune up $r_{\alpha\beta}$by 
recursively,  
so that 
\begin{equation}
e^{ad(\nu \Psi_{\alpha\beta}^*\tilde{h}_{\beta\alpha})}
e^{ad(\nu \Psi_{\alpha\gamma}^*\tilde{h}_{\gamma\beta})}
e^{ad(\nu \Psi_{\alpha\alpha}^*\tilde{h}_{\alpha\gamma})}
=e^{\nu^2c_{\alpha\beta\gamma}^{(2)}ad ( \nu^{-1} ) }~.
\end{equation} 
It follows that $\{\dot\Psi_{\alpha\beta}^*\}$ defines a Weyl manifold 
$\dot{W}_M$ with the Poincar\'{e}-Cartan class 
$c^{(0)}+\nu^2c^{(2)}$. 
Repeating a similar argument as above, we can replace the condition 
mod~$\nu^4$ in (\ref{id3.8}) by mod~$\nu^6$. 
Iterating this procedure, we have 
a Weyl manifold ${W}_M$ such that $c({W}_M)=c
\in \check{H}(M)[[\nu^2]]$. 

Next we would like to show that 
the above construction does not depend on the cocycle chosen. 
Let $\{c_{\alpha\beta\gamma}\},~\{c'_{\alpha\beta\gamma}\} $
be Poincar\'{e}-Cartan cocycles of $\{{C}_U\},
~\{{C}'_U\} $ respectively, which give same 
Poincar\'{e}-Cartan classes. Then, there exists 
$b_{\alpha\beta}\in {\bf R}[[\nu^2]]$ on every $V_\alpha\cap V_\beta\not=
\phi$ such that $b_{\alpha\beta}=-b_{\beta\alpha}$ and 
$c'_{\alpha\beta\gamma}-c_{\alpha\beta\gamma}=b_{\alpha\beta}+b_{\beta\gamma}+b_{\gamma\alpha}.$
Note that 
$b_{\beta\gamma}$ may be replaced by $b_{\beta\gamma}+c_{\beta\gamma}$ 
such that $c_{\alpha\beta}+c_{\beta\gamma}+c_{\gamma\alpha}=0$. 
Since $e^{b_{\alpha\beta}ad(\nu^{-1})}$ is an automorphism, we can replace 
$\Psi_{\alpha\beta}$ by $\acute{\Psi}_{\alpha\beta}=\Psi_{\alpha\beta}
e^{b_{\alpha\beta}ad(\nu^{-1})}$.  
Since $e^{b_{\alpha\beta}ad(\nu^{-1})}$ is the identity on  
${\cal F}({W}_{U_{\alpha\beta}})$, this replacement does not change 
the isomorphism class of ${\cal F}({W}_M)$, but it changes the 
Poincar\'{e}-Cartan cocycle from 
$\{c_{\alpha\beta\gamma}\}$ to $\{c'_{\alpha\beta\gamma}\}. $
This means that the map from the set of Poincar\'{e}-Cartan cocycles into 
the set $\frak W_M$ of isomorphism classes of Weyl manifolds  
induces a map $F$ 
from the set ${\frak{PC}_M}$ of Poincar\'{e}-Cartan classes into $\frak W_M$. 

Next we construct the inverse map  
$\Psi:{W}'_M\rightarrow{W}_M$ which induces the identity on the 
base manifold. 
Equivalently $\Psi^*$ defines an algebra isomorphism of 
${\cal F}({W}_M)$ onto ${\cal F}({W}'_M)$. 
The isomorphism is given by a family $\{\Psi_\alpha^*\}$ of isomorphisms: 
$$\Psi_\alpha^*:{\cal F}({W}_{U_\alpha})
\rightarrow {\cal F}({W}'_{U_\alpha}),$$
each of which induces the identity map on the base space $U_\alpha$ such that 
\begin{equation}\label{id3.5} 
\Psi_\alpha^*\Psi^*_{\alpha\beta}\Psi_\beta^{*-1}=\acute{\Psi}_{\alpha\beta}^*.
\end{equation} 
If we extend $\Psi_\alpha^*$ to a contact Weyl diffeomorphism
$($cf. subsection \ref{subsection-5-3}$)$, then the above 
replacement makes no change of Poincar\'{e}-Cartan cocycle. We use 
the same notation $\Psi_\alpha^*$ for this contact Weyl diffeomorphism. 
By (\ref{id3.5}), and Proposition~\ref{2.4ommyrevised}, we have 
\begin{equation}
\Psi_\alpha^*\Psi^*_{\alpha\beta}\Psi_\beta^{*-1}
e^{b_{\alpha\beta}ad(\nu^{-1})}=\acute{\Psi}_{\alpha\beta}^*.
\end{equation} 
on a contact Weyl algebra bundle. 
However this type of replacement changes the Poincar\'{e}-Cartan cocycle 
within the same cohomology class. 
This means that there is a map from ${\frak W_M}$ 
into ${\frak{PC}_M}$, and which is obviously the inverse map of $F$. 
$\square$
\end{pf}
\par\bigskip 
As mentioned 
in subsection 4.1, there exists a contact Weyl algebra bundle with a 
connection $\nabla^Q$ such that its first Chern class 
coincides with Poincar\'{e}-Cartan class (cf. \cite{y}). 
Then $\nabla^Q|_{{W}_M}$ gives a flat connection on ${W}_M$ and 
there is a one to one correspondence $\sigma$ between 
the space of parallel sections with respect 
to $\nabla|_{{W}_M  } $ and $C^\infty(M)[[\hbar]]$. 
Combining this map with fiberwise star product, we can 
define a star product: 
$f*g=\sigma(\sigma^{-1}(f)*_{\mbox{fiberwise}}\sigma^{-1}(g)). 
$ 
Hence we obtain Theorem~\ref{class}.

 
\section{A Lie group structure of $ {\rm Aut}(M,*)$}
\subsection{Fundamental definitions and representatives}
With the preliminaries in the previous section, 
we give the precise definition of ${\rm Aut}(M,*)$: 
\begin{defn}
\label{MCWD}
\begin{eqnarray}
~~~~~{\rm Aut}(M,*)\!\!\!\!&=&\!\!\!\! 
\{\Psi:C_M\rightarrow C_M\!|\mbox{ fiber-wise }\nu\mbox{-automorphism},
\!\Psi^*({\mathcal F}_M)={\mathcal F}_M\}, \\
~~~~~\underline{\rm Aut}(M,*)\!\!\!\!&=&\!\!\!\!\{\Psi\in {\rm Aut}(M,*)|
   \Psi \mbox{ induces the base identity map. }\} . 
\end{eqnarray}
An element of ${\rm Aut}(M,*)$ is called a {\rm modified 
contact Weyl diffeomorphism}
$($an {\rm MCWD} for short$)$. 
\end{defn}
To illustrate automorphisms of a contact Weyl manifold, 
we consider the automorphisms of a contact Weyl algebra. 
For any real symplectic matrix $A\in Sp(n,{\mathbf R})$,  
set a $\nu$-automorphism of $C$ 
by $\hat{A}Z^i=\sum a_{j}^iZ^j$ and $\hat{A}\nu=\nu$. 
Then we easily have $\hat{A}([a,b])=[\hat{A}a,\hat{A}b]$. 
Conversely, combining the Baker-Campbell-Hausdorff formula 
with the Poincar\'{e} lemma, we have the following. 
\begin{prop}[\cite{y}]\label{y-2-18}
If $\Psi$ is a $\nu$-automorphism of contact Weyl algebra, 
there exists uniquely 
$$\begin{array}{l}
A\in Sp(n,{\mathbf R}),\\ 
F\in \{a=\sum_{2\ell+|\alpha|\ge 3, 
      |\alpha|>0}a_{\ell\alpha}\nu^\ell Z^\alpha\},\\  
c(\nu^2)=\sum_{i=0}^\infty c_{2i}\nu^{2i} \in {\mathbf R}[[\nu^2]],\end{array}$$such that $\Phi=\hat{A}\circ e^{ad(\frac{1}{\nu}(c(\nu^2)+F))}$, 
where $\hat{A}Z^i=\sum a_{j}^iZ^j$ and $\hat{A}\nu=\nu$.
\end{prop}
\noindent{\bf Remark} 
This $\nu$-automorphism can be seen as a ``linear" example appearing  
in the simplest model of contact Weyl manifolds. 
\begin{pf}
The proof of this result was given by Yoshioka \cite{y}, 
however  we recall it because of its importance. 

Let $\frak{m}$ be a unique maximal ideal defined by 
$$\frak{m}=\Bigr\{a=\sum_{2\ell+|\alpha|\ge 1, 
  }a_{\ell\alpha}\nu^\ell Z^\alpha\Bigr\} .$$  
Then we see $\Phi(\frak{m})\subset \frak{m}$, so we have 
\begin{equation}\label{1}
\Phi(Z^i)=\sum a_j^iZ^j+O(2),  
\end{equation}
where $O(2)$ is the collection of the terms degree $\ge 2$. 
Applying (\ref{1}) to the canonical commutation relation 
$[Z^i,Z^j]=\nu\Lambda^{ij}, 
$ 
we see $A\in Sp(n,{\mathbf R})$ and then we may write 
\begin{equation}\label{2}
\hat{A}^{-1}\circ\Phi(Z^i)=Z^i+g_{(2)}^i+O(3),  
\end{equation}
where $g_{(2)}^i$ is the term of homogeneous degree $2$. 
Applying (\ref{2}) to the canonical commutation relation again, 
we have
\begin{eqnarray}
&& \nu\partial_{i+n}g_{(2)}^j=\nu\Lambda^{il}\partial_{l}g_{(2)}^j
=[Z^i,g_{(2)}^j]=[Z^j,g_{(2)}^i]\\
&=&\nonumber
\nu\Lambda^{jk}\partial_{k}g_{(2)}^j
=\nu\partial_{j+n}g_{(2)}^i.
\end{eqnarray}
This is equivalent to 
\begin{equation}
d(\sum_jg_{(2)}^j dz^{j+n}+\sum_i g_{(2)}^idz^{i+n})=0. 
\end{equation}
According to the Poincar\'{e} lemma, there exists uniquely 
$F_{(3)}\in W$ with homogeneous degree 3 such 
that $\frac{1}{\nu}[Z^i,F_{(3)}]=g_{(2)}^i$. 
Therefore we have   
\begin{equation}
\hat{A}^{-1} \circ \Phi (Z^i)= e^{ad(\frac{1}{\nu}F_{(3)})}(Z^i)+O(3), 
\end{equation}
where $O(3)$ is the collection of the terms whose degree $\ge 3$. 
Repeating the above procedure, we obtain that 
\begin{equation}
\hat{A}^{-1} \circ \Phi (Z^i)
= e^{ad(\frac{1}{\nu}F_{(3)})}  
\circ\cdots\circ e^{ad(\frac{1}{\nu}F_{(k)})}(Z^i)+O(k),  
\end{equation}
where $O(k)$ is the collection of the terms whose degree $\ge k$. 
By the Baker-Campbell-Hausudorff formula\footnote{See (\ref{BCH}).}, we have 
\begin{equation}
\hat{A}^{-1} \circ \Phi (Z^i)= e^{ad(\frac{1}{\nu}F)}(Z^i). 
\end{equation}
The uniqueness is inductively 
verified.  
Thanks to the argument above, we may assume that 
$\Psi|_W=\hat{A}\circ e^{ad(\frac{1}{\nu}F)}$. 
Set 
\begin{equation}\label{4}
\tilde{\Phi}(X)=
(\hat{A}\circ e^{ad(\frac{1}{\nu}F)})^{-1} \Psi(X)  .
\end{equation}
Applying (\ref{4}) to $[\tau,Z^i]=\nu Z^i$, 
we have 
\begin{equation}
\tilde{\Psi}(\tau)=\tau+ b(\nu^2)
\end{equation}
for some $b(\nu^2)=b_0+b_2\nu^2+\cdots, \quad b_{2k}\in{\mathbf R}$.
Put $c(\nu^2)=\sum\nu^{2k}\frac{b_{2k}}{2(1-2k)}, $ then we see 
\begin{eqnarray}
&&e^{ad(\frac{1}{\nu}c(\nu^2))}\tau=\tau + b(\nu^2), \\
&&e^{ad(\frac{1}{\nu}c(\nu^2))}Z^i=Z^i.
\end{eqnarray}
Then we have 
\begin{equation}
\Psi=\hat{A}\circ e^{ad(\frac{1}{\nu}(c(\nu^2)+F))}.
\end{equation}
Thus, we see the consequence. 
$\square$\end{pf}

Next we study the basic properties of a modified contact Weyl diffeomorphism. 
We recall the fundamental definitions and facts for readers again. 
Set $\tilde{\tau}_U=\tau+\sum z^i\omega_{ij}Z^j$ 
where $U\subset {\mathbf R}^{2n}$ is an open subset and 
$\omega_{ij}dz^i\land dz^j$ stands for the symplectic structure. 
Then for any modified contact Weyl diffeomorphism, 
we may set ${\Psi}|_{{ C}_U}^\ast(\tilde{\tau}_U)=
a\tilde{\tau}_U+F$, where $a\in C^\infty(U),\, 
F\in\Gamma({ W}_U)$,  
where $W_U$ is a trivial bundle 
$W_U=U\times W$. 
Under the notations above, 
as mentioned before, we have Proposition \ref{representative1}. 
As mentioned in Definition \ref{defn4-10}, 
a modified contact Weyl diffeomorphism ${\Psi}$ 
is called a {\rm contact Weyl diffeomorphism\,(CWD, for short)} if 
\begin{equation}
{\Psi}^\ast\tilde{\tau}_{U'}=\tilde{\tau}_U+f^{\#}. 
\end{equation}
For a contact Weyl diffeomorphism, 
we obtain the following  
(see Corollary 2.5 in \cite{ommy} and Proposition 2.24 in \cite{y}). 
\begin{prop}
\label{2.4ommy}
\begin{enumerate}
\item Suppose 
that ${\Psi}:{C}_U\rightarrow { C}_U$ is a 
contact Weyl diffeomorphism 
which induces the identity map on the base space. 
Then, 
there exists uniquely a Weyl function $f^{\#}(\nu^2)$ of the form 
\begin{equation}
f^{\#}=f_0+\nu^2f_+^{\#}(\nu^2)\quad     
(f_0\in {\bf R},\,f_+(\nu^2)\in C^\infty(U)[[\nu^2]]), 
\end{equation}
such that  
${\Psi} = e^{ad \frac{1}{\nu}\{f_0+\nu^2f_+^{\#}(\nu^2)\}}. $
\item If ${\Psi}$ induces the identity map on $W_U,$ then 
there exists a unique element $c(\nu^2)\in {\bf R}[[\nu^2]]$ with 
$\overline{c(\nu^2)}=c(\nu^2)$, 
such that 
\footnote{Note that this does not induce the 
identity on the whole of ${C}_U$. 
In \cite{y}, a notion of modified contact Weyl diffeomorphism is introduced 
to make a contact Weyl algebra bundle 
$\{{C}_{U_\alpha},
{\Psi}_{\alpha\beta}\}$ 
by adapting the glueing maps to satisfy the cocycle condition and 
patching together them.} 
${\Psi}=e^{ad\frac{1}{\nu}c(\nu^2)} . 
$
\end{enumerate}
\end{prop}
Combining Propositions \ref{representative1} and \ref{2.4ommy}, 
for some $g(\nu^2)$, we see that 
$$
\Psi\circ 
e^{ad(\frac{1}{\nu}g(\nu^2))}=e^{ad(\frac{1}{\nu}f^\#(\nu^2))},     
$$
where 
$f^{\#}(\nu)$ is a Weyl function with the form 
$$
f^{\#}(\nu^2)=f_0+\nu^2f_+^{\#}(\nu^2)\quad     
(f_0\in {\bf R},\,f_+(\nu^2)\in C^\infty(U)[[\nu]]).    
$$
Thus we obtain 
$$
\Psi=e^{ad(\frac{1}{\nu}(g(\nu^2)+f^\#(\nu^2)))}.
$$
The following is easily verified. 
\begin{prop}\label{lemma1}
Suppose that $f(\nu^2),~ a(\nu^2)\in C^\infty(U)[[\nu^2]]$. 
If 
\begin{equation}
\label{ass}
\nu^2f^{\#}(\nu)=g(\nu^2),
\end{equation}
then we have 
$g_0=0$ and $f(\nu^2)=\sum_{i\ge0}g_i\nu^{2i} \in {\bf R}[[\nu^2]]\nu^2$. 
\end{prop}
\begin{pf}
The right hand side of (\ref{ass}) has no term containing a 
factor $X^\alpha Y^\beta\,(|\alpha|+|\beta|\ge1)$.   
Hence $\partial_X^\alpha\partial_Y^\beta f =0\, (|\alpha|+|\beta|\ge1)$. 
Then we see 
$f=\sum c_i\nu^{2i}$. Thus, 
$\nu^2\sum c_i\nu^{2i}=\sum g_i\nu^{2i}$ and $g_0=0$. 
$\square$\end{pf}
Using the above proposition, we have 
\begin{prop}\label{lemma3}
If 
\begin{equation}\label{ass3}
e^{ad(\frac{1}{\nu}\{g+\nu^2 f^{\#}\})}=e^{ad(\frac{1}{\nu}\{g'+\nu^2
 f^{'\#}\})} , \end{equation}
then we have  
\begin{equation}\label{conc3}
   \{g+\nu^2 f^{\#}\} = \{g'+\nu^2f^{'\#}\} 
\end{equation}
and 
\begin{equation}
\nu^2(f-f')=g'-g\in \nu^2{\bf R}[[\nu^2]]  . 
\end{equation}\end{prop}  
\begin{pf}
Applying $e^{ad(\frac{1}{\nu}g)}$ 
to the both hand side, we have 
\begin{equation}
e^{ad(\frac{1}{\nu}\{\nu f^{\#}\})}=e^{ad(\frac{1}{\nu}\{-g+g'+\nu f^{'\#}\})},
\end{equation}
and this implies 
that the left hand side of the above equality is a contact Weyl diffeomorphism. Then by uniqueness in Proposition\,\ref{2.4ommy} and Proposition\,
\ref{lemma1}, 
we see the consequence. 
$\square$\end{pf} 
We have 
\begin{prop}\label{representative2}
For any modified contact Weyl diffeomorphism 
$\Psi:C_U \rightarrow C_{U}$ which induces the identity map on the base space, 
there exists a Weyl function $f^{\#}(\nu^2)$ of the form 
\begin{equation}
f^{\#}(\nu^2)=f_0+\nu^2f_+^{\#}(\nu^2)\quad     
(f_0\in {\bf R},\,f_+(\nu^2)\in C^\infty(U)[[\nu^2]]), 
\end{equation}
and smooth function $g(\nu^2)\in C^\infty(U)[[\nu^2]]$ 
such that 
$\Psi=e^{ad(\frac{1}{\nu}\{g(\nu^2)+f^\#(\nu^2)\} ) } .$
\end{prop}
\noindent{\bf Remark} 
Please compare this result with Proposition~\ref{y-2-18}. 
\par\medskip \noindent
Furthermore, we have 
\begin{prop}\label{representative} 
Let $\Psi_{U_\alpha}$ {\rm (}resp. $\Psi_{U_\beta}${\rm )} be 
a modified contact 
Weyl diffeomorphism on $C_{U_\alpha}$ {\rm (}resp. $C_{U_\beta}${\rm )} 
inducing the identity map on the base manifold. 
Suppose that  
$${\Psi}_{U_\alpha}|_{C_{U_{\alpha \beta}}}=
{\Psi}_{U_\beta}|_{C_{U_{\beta \alpha}}},   
$$
where 
$U_{\alpha\beta}:=\varphi_\alpha(V_\alpha\cap V_\beta), 
~U_{\beta\alpha}:=\varphi_\beta(V_\alpha\cap V_\beta),~ 
C_{U_{\alpha\beta}}:=\Psi_\alpha(C_{V_\alpha}|_{V_\alpha\cap V_\beta})$ 
etc\footnote{See also Definition \ref{cwm}.}. 
Then 
\begin{equation}
\Psi_{\alpha\beta}^{-1,\ast}
\bigr((g_\alpha(\nu^2)+\nu^2 f_\alpha^{\#}(\nu^2))|_{U_{\alpha\beta}}\bigr)=
(g_\beta(\nu^2)+\nu^2 f_\beta^{\#}(\nu^2))|_{U_{\alpha\beta}}.
\end{equation}
Thus, patching $\{g_U+\nu^2f_U^{\#}\}$ together  
we can make a global function 
$g+\nu^2f^{\#} \in 
C^\infty(M)[[\nu^2]]+\nu^2C^\infty(M)^\#[[\nu^2]]
$. 
Hence there is a bijection between 
$\underline{\rm Aut}(M,*)$ 
and $C^\infty(M)[[\nu^2]]+\nu^2C^\infty(M)^{\#}[[\nu^2]]$.
\end{prop}
The propostions mentioned above implicate that 
the space 
$$\frak{C}_c(M)=C_c^\infty(M)[[\nu^2]]+\nu^2C_c^\infty(M)^\#[[\nu^2]]$$ 
is a candidate of the model space of $\underline{\rm Aut}(M,*)$. 
In fact, 
the Baker-Campbell-Hausdorff formula shows the smoothness of 
group operations.   
Therefore we have the following: 
\begin{thm}\label{Lie-group-underline-aut}
$\underline{\rm Aut}(M,*)$ is a Lie group 
modeled on 
$\frak{C}_c(M)$. 
\end{thm}
\begin{pf}
The smoothness of group operations is ensured by the following formula:  
Put $H_i(\nu^2)=g_i(\nu^2)+\nu^2 f_i^\#(\nu^2)$ $(i=1,2).$ 
\begin{eqnarray}\label{BCH}
&&e^{ad(\frac{1}{\nu}H_1(\nu^2))}
\circ e^{ad(\frac{1}{\nu}H_2(\nu^2))}\\ \nonumber
&=& \!\!\!\! e^{ad[\frac{1}{\nu}\{(H_1(\nu^2)+H_2(\nu^2))+\cdots+
\frak{B}_m(H_1(\nu^2),H_2(\nu^2))+\cdots \}]
},
\end{eqnarray}
where $\frak{B}_m$ means 
the general term of the Baker-Campbell-Hausdorff formula
\footnote{More 
precisely $\frak{B}_m$ is given by the following way: 
\\
$\quad\qquad\qquad
\frak{B}_m(\frac{1}{\nu}H_1(\nu^2),\frac{1}{\nu}H_2(\nu^2))
$
\\
$\qquad\qquad
\!\!=\!\!\!
\scriptsize {\frac{(-1)^{m-1}}{m}\sum\frac{ad(\frac{1}{\nu}H_1(\nu^2))^{p_1}
ad(\frac{1}{\nu}H_2(\nu^2))^{q_1}
\cdots ad(\frac{1}{\nu}H_1(\nu^2))^{p_m}ad(\frac{1}{\nu}
H_2(\nu^2))^{q_m-1}\frac{1}{\nu}H_2(\nu^2)}
{ p_1 !q_1!\cdots p_m!q_m!( p_1+ q_1+ \cdots + p_m+ q_m )}.
}
$}.
\begin{eqnarray}
\left\{e^{ad(\frac{1}{\nu}H_1(\nu^2))}\right\}^{-1}=
e^{ad(\frac{-1}{\nu}H_1(\nu^2))} .
\end{eqnarray}
This completes the proof. $\square$\end{pf}
%
As to $\frak{C}_{c}(M)$, we have 
\begin{lem}
\label{star-quasi-multiplicative}
The space 
$$\frak{C}_c(M)={\rm ind}\lim_{K\mbox{:compact}}
\Bigl(C_K^\infty(M)[[\nu^2]]+\nu^2C_K^\infty(M)^\#[[\nu^2]]\Bigl)
$$ is Mackey complete, 
where $C_K^\infty(M)$ is the space equipped with the standard 
locally convex topology. 
\end{lem}
\begin{pf}
Remark that $C^\infty(M,N)$ is a smooth manifold modeled on a Mackey complete 
locally convex space 
$C_c^\infty(M\leftarrow f^*TN)$, where $M$ and $N$    
are finite-dimensional manifolds. 
Since limits, direct sums and inductive limits preserve 
Mackey completeness, $\frak{C}_{c}(M)$ is also a 
Mackey complete 
locally convex space.
$\square$\end{pf}
In general, we can show the followings. 
\begin{lem}\label{lemma-quasi-multiplicative}
If $(E,~\star,~||\cdot||_\rho)$ is a 
Mackey complete locally convex space with a quasi multiplicative\footnote{
The assumption (\ref{quasi-multiplicative}) can be replaced by 
$||f\star g||_\rho \le C_\rho ||f||_\rho\cdot||g||_{\hat{\rho}} $.}
binary operation $\star$, that is, 
\begin{equation}
\label{quasi-multiplicative}
||f\star g||_\rho \le C_\rho ||f||_\rho\cdot||g||_\rho, 
\end{equation}
for some positive number $C_\rho$, 
then $\sum_{n=0}^\infty\frac{f\star\cdots\star f}{n!}$ converges. 
Set $e_\star^{f}=\sum_{n=0}^\infty\frac{f\star\cdots\star f}{n!}$  
Then we have 
\begin{equation}\label{star-inequality}
||e_\star^f ||_\rho \le \sum \frac{C_\rho^{n-1}||f||\rho^n}{n!}. 
\end{equation}
\end{lem}
\begin{pf}
By the assumtion, we have 
\begin{eqnarray}
||f\star\cdots\star f||_\rho&\le& 
C_\rho ||f||_\rho \cdots ||f \star \cdots \star f||_\rho \\
\nonumber
&\le& 
C_\rho^{n-2}||f||_\rho^{n-2} ||f \star f||_\rho \\ \nonumber
&\le& 
C_\rho^{n-1}||f||_\rho^n.
\end{eqnarray}
Hence we see that 
$\left\{
\sum_{n=0}^\ell \frac{f \star \cdots \star f}{n!}\right
\}_{n=1}^\infty$ 
is a Mackey-Cauchy sequence. 
Set 
$
B=\left\{\sum_{n=0}^\ell \frac{f \star \cdots \star f}{n!}
\right\}_{\ell=0}^\infty  .$  
Then by the Mackey completeness, $E_B$ is complete. Hence there 
exists uniquely an element denoted by 
$e_\star^{f}=
\sum_{n=0}^\infty \frac{f \star \cdots \star f}{n!} $ 
such that 
$$\sum_{n=0}^\ell \frac{f \star \cdots \star f}{n!} 
\rightarrow e_\star^{f} \in E_B \subset E , $$ 
and we also have 
$$||\sum_{n=0}^\infty \frac{f \star \cdots \star f}{n!} ||_\rho 
\le \sum_{n=0}^\infty \frac{C_\rho^{n-1}||f||_\rho^n}{n!} .
$$ 
$\square$\end{pf}
\begin{lem}\label{lemma-product-integral}
Let 
$(E,~\star,~||\cdot||_\rho)$ be a Mackey complete 
locally convex space  
with a quasi multiplicative binary operation. 
Then for any smooth curve $X(t)$ in $E$, the product integral 
\begin{equation}\label{product-integral}
\prod e_\star^{X(t)dt}=\lim_{n\rightarrow\infty }
                      e_\star^{X(t_n)\Delta t_n}\star\cdots\star 
                      e_\star^{X(t_i)\Delta t_i}\star\cdots\star
                      e_\star^{X(t_0)\Delta t_0}
\end{equation}
exists. 
\end{lem}
\begin{pf}
For $I=[t_0,t_n]$, set 
$\Delta:t_0<\cdots < t_i < \cdots t_n$, $\Delta t_i:=t_{i+1}-t_i$, and 
${\rm mesh}(\Delta) : =\max\{\Delta t_0,\ldots,
\Delta t_i,\ldots,\Delta t_{n-1}\}.$ 
We have to show 
$$      e_\star^{X(t_n)\Delta t_n}\star\cdots\star 
                      e_\star^{X(t_i)\Delta t_i}\star\cdots\star
                      e_\star^{X(t_0)\Delta t_0}
$$ is a Mackey-Cauchy net. 
A direct computation gives the following estimation: 
\begin{eqnarray}
\label{keyinequality0} 
&&~~~\Bigr|\Bigr|e_\star^{X(t_{n-1})\Delta t_{n-1}}\star\cdots\star 
                      e_\star^{X(t_i)\Delta t_i}\star\cdots\star
                      e_\star^{X(t_0)\Delta t_0}
\\\nonumber
&&\qquad -e_\star^{X(t_0)\Delta t_{n-1}}\star\cdots\star 
                      e_\star^{X(t_0)\Delta t_i}\star\cdots\star
                      e_\star^{X(t_0)\Delta t_0}\Bigr|\Bigr|_\rho \\
&=& 
\nonumber 
\Bigr|\Bigr|\sum_{i=0}^{n-1}
\left(e_\star^{X(t_{0})\Delta t_{n-1}}\star\cdots\star 
                      e_\star^{X(t_{0})\Delta t_{i+1}}
                      e_\star^{X(t_0)\Delta t_i}
                      e_\star^{X(t_{i-1})\Delta t_{i-1}}\star\cdots\star
                      e_\star^{X(t_0)\Delta t_0}\right.\\
&&\nonumber \qquad \left.-
e_\star^{X(t_{0})\Delta t_{n-1}}\star\cdots\star 
                      e_\star^{X(t_{0})\Delta t_{i+1}}
                      e_\star^{X(t_i)\Delta t_i}
                      e_\star^{X(t_{i-1})\Delta t_{i-1}}\star\cdots\star
                      e_\star^{X(t_0)\Delta t_0}
\right)
\Bigr|\Bigr|_\rho
\\
&\le&\nonumber
\Bigr|\Bigr|\sum_{i=0}^{n-1}
\left(e_\star^{X(t_{0})\Delta t_{n-1}}\star\cdots\star 
                      e_\star^{X(t_{0})\Delta t_{i+1}}\star\right.\\
&&\nonumber
                      \left.
                   (e_\star^{X(t_0)\Delta t_i}-e_\star^{X(t_i)\Delta t_i})\star
                   e_\star^{X(t_{i-1})\Delta t_{i-1}}\star\cdots
                   \star
                   e_\star^{X(t_0)\Delta t_0}\right)\Bigr|\Bigr|_\rho\\
&\le&\nonumber
\sum_{k=0}^{n-1}
e^{\sup_{t\in I=[t_0,t_n]} || X(t)||_\rho \Delta t_{n-1}}
\times\cdots\times
\\
&&\nonumber\qquad
\cdots \times
\Bigr|\Bigr|e^{\sup_{t\in I=[t_0,t_n]} X(t_0)\Delta t_{k}}
 -e^{\sup_{t\in I=[t_0,t_n]} X(t_k)\Delta t_{k}}\Bigr|\Bigr|_\rho
\times\cdots
\\
&&\nonumber\qquad\qquad\qquad\times\cdots\times
e^{\sup_{t\in I=[t_0,t_n]} || X(t)||\rho \Delta t_{0}}
\\
&\le& \nonumber
\sum_{k=0}^{n-1}
\Bigr|\Bigr|X(t_k)-X(t_0)\Bigr|\Bigr|_\rho
|\Delta t_k|e^{\sup_{t\in I}||X(t)||_\rho}
\\ 
&\le& \nonumber 
|I|^2\sup_{t\in I}||X'(t)||_\rho e^{|I|\sup_{t\in I}||X(t)||_\rho}
.
\end{eqnarray}
in the last inequality, we used (\ref{star-inequality}).
Let $\Delta:a=s_0< \cdots < s_\ell < \cdots< s_m=b$ be a  
division of $[a,b]$ and 
$\Delta^{(\ell)}
:s_\ell=
t_0^{(\ell)}<\cdots < t_i^{(\ell)}<\cdots t_{n_{(\ell)}}^{(\ell)}=s_{\ell+1}$ 
a subdivision of $\Delta$. 
Then 
\begin{eqnarray}
\nonumber
&&\Bigr|\Bigr|e_\star^{X(s_{m-1})\Delta s_{m-1}}\star\cdots\star 
                      e_\star^{X(s_i)\Delta s_i}\star\cdots\star
                      e_\star^{X(s_0)\Delta s_0}\\
&&\nonumber \quad
-\left( e_\star^{X(t_{n_{m-1}-1}^{(m-1)})\Delta t_{n_{m-1}-1}^{(m-1)}}\star
\cdots\star 
                e_\star^{X(t_i^{(m-1)})\Delta t_i^{(m-1)}}\star\cdots\star
                e_\star^{X(t_0^{(m-1)})\Delta t_0^{(m-1)}}\right)\star\cdots\star\\
   &&\nonumber
   \qquad\quad\star\cdots\star
 \left(e_\star^{X(t_{n_{\ell}-1}^{(\ell)})\Delta t_{n_{\ell}-1}^{(\ell)}}\star
 \cdots\star                 
 e_\star^{X(t_i^{(\ell)})\Delta t_i^{(\ell)}}\star\cdots\star
                e_\star^{X(t_0^{(\ell)})\Delta t_0^{(\ell)}}\right)\star\cdots
                \star\\
   &&\nonumber
   \qquad\qquad\qquad\star\cdots\star             
   \left(e_\star^{X(t_{n_{0}-1}^{(0)})\Delta t_{n_{0}-1}^{(0)}}\star\cdots
   \star 
                e_\star^{X(t_i)^{(0)}\Delta t_i^{(0)}}\star\cdots\star
                e_\star^{X(t_0)^{(0)}\Delta t_0^{(0)}}\right)\Bigr|\Bigr|_\rho
\\
\nonumber
&\le&
\sum_{\ell=0}^{m-1}\Bigr|\Bigr|
e_\star^{X(s_{m-1})\Delta s_{m-1}}\star\cdots\star \\
&&\nonumber\quad\qquad\star\cdots\star
\left( e_\star^{X(s_{\ell})\Delta s_{\ell}}
      -\left(e_\star^{X(t_{n_{\ell}-1}^{(\ell)})\Delta t_{n_{\ell}-1}^{(\ell)}}\star\cdots\star                 
      e_\star^{X(t_i^{(\ell)})\Delta t_i^{(\ell)}}\star\cdots\star
e_\star^{X(t_0^{(\ell)})\Delta t_0^{(\ell)}}\right)\right)\star\cdots\star\\
&&\nonumber\qquad\qquad\quad\star\cdots\star 
    \left(e_\star^{X(t_{n_{0}-1}^{(0)})
                \Delta t_{n_{0}-1}^{(0)}}\star\cdots\star 
                 e_\star^{X(t_i^{(0)})\Delta t_i^{(0)}}\star\cdots\star
                e_\star^{X(t_0^{(0)})\Delta t_0^{(0)}}\right)\Bigr|\Bigr|_\rho  \\   
 &\stackrel{(*)}{\le}&
\nonumber
\sum_{\ell=0}^{m-1}e^{(\sup_{s\in [a,b]} ||X(s)||_\rho)\Delta s_{m-1}}
\cdots \\
&&\nonumber\qquad\quad\qquad
\cdots |\Delta s_\ell|^2
{(\sup_{s\in [a,b]} ||X'(s)||_\rho)}
e^{(\sup_{s\in [a,b]} ||X(s)||_\rho)\Delta s_{\ell}}
\cdots e^{(\sup_{s\in [a,b]} ||X(s)||_\rho)\Delta s_{0}}\\
&\le&\nonumber
(\sup_{s\in [a,b]} ||X'(s)||_\rho)(b-a)
\max_{\ell}|\Delta s_\ell|e^{(\sup_{s\in [a,b]} ||X(s)||_\rho)(b-a)}.
\end{eqnarray}
In the estimation $(*)$, we used (\ref{keyinequality0}).  
This implies that 
$$  \{    e_\star^{X(s_m)\Delta s_m}\star\cdots\star 
                      e_\star^{X(s_\ell)\Delta s_\ell}\star\cdots\star
                      e_\star^{X(s_0)\Delta s_0}   \}_\Delta
$$ is a Mackey-Cauchy net. 
$\square$\end{pf}

Before stating the next lemma, 
we recall the precise definition of seminorms which we use. 
The seminorms of $C_c^\infty(M)^\#[[\nu^2]]$ are defined by: 
\begin{equation}\label{seminorm}
\bigr|\bigr|\sum_{l=0}^\infty\nu^\ell f_\ell \bigl|\bigl|_{i,K}:=
\sum_{|\alpha|+ 2\ell \le i}
\sum_{p\in K}|\partial_z^\alpha f_\ell(p)|, ~~~(i\in {\mathbf N})
\end{equation}
where $K$ is a compact subset of $M.$
Then we have the following. 
\begin{lem}\label{star-product-quasi-multiplicative}
Set $f(\nu)=\sum_{k\in {\mathbf N}}f_k\nu^k$ and 
$g(\nu)=\sum_{\ell\in {\mathbf N}}f_\ell \nu^\ell$. 
\begin{equation}
\bigr|\bigr|f(\nu)*g(\nu)\bigl|\bigl|_{i,K}\le 
C_{i,K}\bigr|\bigr|f(\nu) \bigl|\bigl|_{i,K}\bigr|\bigr|
g(\nu) \bigl|\bigl|_{i,K}.
\end{equation}
\end{lem}
\begin{pf}
We may assume that $K$ is a subset of a Darboux chart $(U;(x,\xi))$. 
\begin{eqnarray}\label{Moyal-concrete}
&&
\bigr|\bigr|f(\nu)* g(\nu)\bigl|\bigl|_{i,K} \\
&=&\nonumber \Bigr|\Bigr|f(\nu)
e^{\frac{\nu}{2}\stackrel{\leftarrow}{\partial_x}\land 
\stackrel{\rightarrow}{\partial_\xi}}g(\nu) 
\Bigr|\Bigr|_{i,K}
\\
\nonumber
&=&\sum_{\alpha\beta}  
\Bigr|\Bigr|
\frac{\left(\frac{\nu}{2}\right)^{k+\ell+|\alpha+\beta|}}{\alpha!\beta!}
  \partial_x^\alpha\partial_\xi^\beta f\cdot
     \partial_x^\beta(-\partial_\xi)^\alpha g \Bigr|\Bigr|_{i,K}\\
\nonumber
&=&
\sum_{2k+2\ell + 2|\alpha|+2|\beta|+|\gamma|+|\delta|\le i}
\sup_K
\left|
\frac{\left(\frac{1}{2}\right)^{|\alpha+\beta|}}{\alpha!\beta!}
\partial_x^\gamma\partial_\xi^\delta
(  \partial_x^\alpha\partial_\xi^\beta f_k\cdot
     \partial_x^\beta(-\partial_\xi)^\alpha g_\ell )\right|\\
\nonumber
&=&\sum_{2k+2\ell + 2|\alpha|+2|\beta|+|\gamma|+|\delta|\le i}
\Bigr|\frac{\left(\frac{1}{2}\right)^{|\alpha+\beta|}}{\alpha!\beta!}
\sum_{\zeta,\eta}{}_\gamma C_\zeta {}_\delta C_\eta  
\\
\nonumber& &\qquad\qquad\qquad\qquad \times      
\bigr(  \partial_x^{\alpha+\zeta}\partial_\xi^{\beta+\eta} f_k\cdot
\partial_x^{\beta+(\gamma-\zeta)}(-\partial_\xi)^{\alpha+(\delta-\eta)}g_\ell
\bigr)\Bigl|
\\
\nonumber 
&\le&
\sum_{2k+2\ell + 2|\alpha|+2|\beta|+|\gamma|+|\delta|\le i}
\frac{\left(\frac{1}{2}\right)^{|\alpha+\beta|}}{\alpha!\beta!}
\sum_{\zeta,\eta}{}_\gamma C_\zeta {}_\delta C_\eta
||f ||_{i,K}\cdot||g ||_{i,K}.  
\end{eqnarray}
In the last estimation, we used the followings. 
\begin{eqnarray}
\nonumber
|\alpha|+|\zeta|+|\beta|+|\eta|+k 
\nonumber
&\le&
2|\alpha|+|\gamma|+2|\beta|+|\delta|+2k+2\ell, 
\end{eqnarray}
\begin{eqnarray}
\nonumber
|\alpha|+|\gamma-\zeta|+|\beta|+|\delta-\eta|+\ell  
\nonumber
&\le&
2|\alpha|+|\gamma|+2|\beta|+|\delta|+2k+2\ell.
\end{eqnarray}
$\square$\end{pf}
Using this lemma, we easily have 
\begin{lem}\label{star-quasi-multiplicative-cor}
For any elements, 
$H_i(\nu^2)=g_i(\nu^2)+\nu^2 f_i^\#(\nu^2) \in \frak{C}_{c}(M)
=C^\infty_c(M)[[\nu^2]] + \nu^2C_c^\infty(M)^\#[[\nu^2]]~(i=1,2),$ 
define a product in the following way.
\begin{eqnarray}\label{BCH_defn}
&& (\frac{1}{\nu} H_1(\nu^2)) \star (\frac{1}{\nu} H_2(\nu^2))\\
&=& \frac{1}{\nu}\{(H_1(\nu^2)+H_2(\nu^2))+\cdots+
\frak{B}_m(H_1(\nu^2),H_2(\nu^2))+\cdots \},
\nonumber
\end{eqnarray}
where $\frak{B}_m$ denotes 
the general term of the Baker-Campbell-Hausdorff formula.
Then ${\rm Aut}(*)\cong (\frak{C}_{c}(M),\star)$ as a 
Mackey complete locally convex space with a quasi multiplicative 
binary operation. 
\end{lem}

Making use of Lemmas~\ref{star-quasi-multiplicative}, 
\ref{lemma-quasi-multiplicative} and \ref{star-quasi-multiplicative-cor},  
we can show the exsistence  
of solution for the equation (\ref{regularity-eq}) 
when $G=\underline{\rm Aut}(M,\ast)$ in Definition~\ref{regular} 
( cf. \cite{tokyo-miya} and \cite{miy}). 
Then we see that smooth 
curves in the Lie algebra 
integrate to smooth curves in the Lie group in a smooth way. 
Thus we have  
\begin{thm}\label{regularity-underline-aut}
$\underline{\rm Aut}(M,*)$ is a regular Lie group 
modeled on 
$\frak{C}_c(M)$. 
\end{thm}

\subsection{Lifts as modified contact Weyl diffeomorphisms}
As will be seen in the next proposition, 
general modified contact Weyl diffeomorphims are closely related to 
symplectic diffeomorphisms. 
\begin{prop}
\label{inducing}
For any modified contact Weyl diffeomorphism $\Psi$, it induces 
a symplectic diffeomorphism 
on the base symplectic manifold. 
Moreover, there exists a group homomorphism 
$p$ from ${\rm Aut}(M,*)$ into ${\rm Diff}(M,\omega)$. 
\end{prop} 
Conversely, we consider the following problem: 
\par\bigskip\noindent
{\bf Problem} 
{\it For any} globally {\it defined symplectic diffeomorphism 
$\phi:M\rightarrow M$, does there exist a} globally {\it defined 
modified contact Weyl  
diffeomorphism~{\rm (}referred to as a {\it MCW-lift}{\rm )} 
$\hat{\phi}$ which induces $\phi$~?}. 
\par\bigskip\noindent
To solve the problem above, we need several notations. 
Let $(M,\omega)$ be a symplectic manifold and $W_M$ the Weyl algebra bundle 
over $(M,\omega)$. 
Set 
\begin{eqnarray}
&&\nabla^{symp}:=\mbox{the canonical extention of symplectic connection to }
W_M , \\
&&\delta:=ad(\frac{1}{\nu}\omega_{ij}dz^iZ^j) , \\
&&\nabla^W:=\nabla|_{W_M}=\nabla^{symp}-\delta+ad(\frac{1}{\nu}\gamma):
\mbox{a Fedosov connection} , \\ 
&&\label{phi}
\phi:M\rightarrow M:\mbox{ a symplectic diffeomorphism},\\
&&\label{G}
\tilde{\nabla}^W:=\nabla^{symp}+ad(\frac{1}{\nu}\phi^{-1*}G),\quad
(\mbox{ where }G:=\omega_{ij}dz^i Z^j+\gamma) , 
\\
&&\{i;D;j\}\bigl(\frac{F}{\nu}\bigl):=\Bigr(
\bigl(ad(\frac{F}{\nu})\bigl)^i\circ ad(D\bigl(\frac{F}{\nu}\bigl)) \circ
\bigl(ad(\frac{F}{\nu})\bigl)^j\Bigl). \label{i-D-j}
\end{eqnarray}
Here we remark that for any symplectic diffeomorphism $\phi$ on $M$ and any 
section $\sigma\in W_M$, the pull-back $\phi^*(\sigma)$ is naturally extended 
in the following way. 
\begin{equation}
\phi^*(\sigma(z,\nu,dz))=\sigma(\phi^*(z),\nu,\phi^*(dz)) . 
\end{equation}
In order to construct a lift of a symplectic diffeomorphism, we need 
several formulas.  
\begin{lem}
Under the notations above,
\begin{equation}
\nabla^W\circ \phi^*(\sigma(z,\nu,dz))
=\phi^*\circ(\nabla^W)
-\phi^*\circ\Bigr(ad(\frac{1}{\nu}G_z-\frac{1}{\nu}\phi^{-1*}(G)|_z)\Bigr).
\label{star}
\end{equation}
\end{lem}
\begin{pf}
By a direct computation, we have 
\begin{eqnarray}
&&\nabla^W\circ \phi^*(\sigma(z,\nu,dz))
\\
&=&\nonumber
\Bigr(\bigr\{\nabla^{symp}+ad(\frac{1}{\nu}G_y)\bigr\} \circ \phi^*\Bigr)
\sigma(z,\nu,dz)
\\
&=&\nonumber
\bigr(\nabla^{symp}\circ\phi^*)\sigma(z,\nu,dz)\\
&&\nonumber\qquad\qquad +ad(\frac{1}{\nu}G_y)\phi^*(\sigma)(y,\nu,dy)\\
&\stackrel{\mbox{symp. conn.}}{=}&\nonumber 
\bigr(\phi^*\circ \nabla^{symp})\sigma(z,\nu,dz)\\
&&\nonumber\qquad\qquad +ad(\frac{1}{\nu}G_y)\phi^*(\sigma)(y,\nu,dy)\\
&=&\nonumber
\bigr(\phi^*\circ\nabla^{symp})\sigma(z,\nu,dz)\\
&&\nonumber\qquad\qquad + 
\phi^*\Bigr( ad(\frac{1}{\nu}(\phi^{-1*}G_y)|_z ) \sigma(z,\nu,dz) \Bigr)\\
&=&\nonumber
(\phi^*\circ(\nabla^{symp}+ad(\frac{1}{\nu}G_z) ) )\sigma(z,\nu,dz)\\
&&\nonumber\qquad 
-\phi^*\circ ad(\frac{1}{\nu}G_z)\sigma(z,\nu,dz)\\
&&\nonumber\qquad\qquad
+ 
\phi^*\Bigr( ad(\frac{1}{\nu}(\phi^{-1*}G_y)|_z ) \sigma(z,\nu,dz) \Bigr)
\\%
&=&\nonumber(\phi^*\circ\nabla^W )\sigma(z,\nu,dz)\\
\nonumber
&&\qquad
-\phi^*\circ \Bigr( ad(\frac{1}{\nu}G_z)\sigma(z,\nu,dz)
- ad(\frac{1}{\nu}(\phi^{-1*}G_y)|_z ) \Bigr) \sigma(z,\nu,dz). 
\end{eqnarray}
Then we obtain the desired formula. $\square$
\end{pf}
\noindent
According to the lemma above, we have 
\begin{lem}
Under the same notations above, 
for $n > 0,$ 
\begin{eqnarray}
\nonumber 
\nabla^W\circ(ad\frac{F}{\nu})^n 
=\nonumber\sum_{i+j=n-1}(ad\frac{F}{\nu})^i\circ
ad(\nabla^W\frac{F}{\nu})\circ 
ad(\frac{F}{\nu})^j 
+ 
(ad\frac{F}{\nu})^n\circ \nabla^W. 
\end{eqnarray}
\end{lem}
\begin{pf}
First, we can easily verify that  
\begin{eqnarray}
&&\nonumber 
\nabla^W\circ(ad\frac{F}{\nu})\sigma 
=\nonumber\nabla^W[\frac{F}{\nu},\sigma]
\stackrel{alg. conn.}{=}\nonumber[\nabla^W \frac{F}{\nu},\sigma]
+[\frac{F}{\nu},\nabla^W\sigma],  
\end{eqnarray}
for any section $\sigma \in \Gamma(W_M)$. 
By using induction on $n$, we see that 
\begin{eqnarray}
\nonumber 
\nabla^W\circ(ad\frac{F}{\nu})^n\sigma 
=\nonumber\sum_{i+j=k-1}(ad\frac{F}{\nu})^i\circ
ad(\nabla^W\frac{F}{\nu})\circ 
ad(\frac{F}{\nu})^j \sigma
+ 
(ad\frac{F}{\nu})^k\circ \nabla^W\sigma. 
\end{eqnarray}
This completes the proof. $\square$
\end{pf}
\noindent
We also have 
\begin{lem}
Set $\{i;\nabla^W;j\}=(ad\frac{F}{\nu})^i\circ 
ad(\nabla^W\frac{F}{\nu})\circ(ad\frac{F}{\nu})^j$ following (\ref{i-D-j}). 
Then 
\begin{equation}~~~
\nabla^W\circ \exp[ad(\frac{F}{\nu})]
\!=\!\exp[ad(\frac{F}{\nu})]\circ\nabla^W \!\!+\!
\sum_{k=1}^{\infty}\frac{1}{k!}
\Bigr(\sum_{i+j=k-1}\{i;\nabla^W;j\}(\frac{F}{\nu})\Bigr).   
\label{star-star}
\end{equation}
\end{lem}
\begin{pf}
By a direct computation, we see that 
\begin{eqnarray}
&&\nonumber\nabla^W\circ\exp[ad\frac{F}{\nu}]\\
&=&\nonumber
\nabla^W\circ \sum_{k=0}^{\infty} \frac{1}{k!}(ad\frac{F}{\nu})^k\\
&=&\nonumber
\sum_{k=1}^{\infty} \frac{1}{k!}
\Bigr(\sum_{i+j=k-1}(ad\frac{F}{\nu})^i\circ
(ad(\nabla^W\frac{F}{\nu}))\circ(ad\frac{F}{\nu})^j\Bigr)
+(ad\frac{F}{\nu})^k\circ\nabla^W .
\end{eqnarray}
Thus we obtain the desired formula. $\square$
\end{pf}
Thanks to the lemmas above, we obtain the following.  
\begin{thm}\label{lifting}
Under the same notations, 
\begin{eqnarray}
\nonumber
\nabla^W\circ\phi^*\circ\exp[ad(\frac{F}{\nu})] 
&=&\nonumber
\phi^*\circ\exp[ad(\frac{F}{\nu})] \circ\nabla^W 
+ \phi^*\sum_{k=1}^{\infty} \sum_{i+j=k-1}\{i;\tilde{\nabla}^W;j\}
(\frac{F}{\nu})\\
&&\nonumber
\qquad\qquad-\phi^*\Bigr(\exp[ad(\frac{F}{\nu})]\bigr(\frac{1}{\nu}(
G_z-\phi^{-1*}(G_y)_z)\bigr)\Bigr). 
\end{eqnarray} 
\end{thm} 
\begin{pf} 
A direct computation with formulas ~(\ref{star}) and (\ref{star-star}) 
gives 
\begin{eqnarray}
&&\nonumber
\nabla^W\circ\phi^*\circ\exp[ad(\frac{F}{\nu})] \\
&\stackrel{(\ref{star})}{=}&\nonumber
\bigr\{
\phi^*\circ \nabla^W - \phi^* \circ 
ad(\frac{1}{\nu}(G_z-\phi^{-1*}(G_y)_z) )  
\bigr\} 
\circ\exp[ad(\frac{F}{\nu})] \\ 
&=&\nonumber
\phi^*\circ \nabla^W \circ\exp[ad(\frac{F}{\nu})] 
- \phi^* \circ 
ad(\frac{1}{\nu}(G_z-\phi^{-1*}(G_y)_z) )  \circ\exp[ad(\frac{F}{\nu})] \\
&\stackrel{(\ref{star-star})}{=}&\nonumber
\phi^*\circ\exp[ad(\frac{F}{\nu})]\circ\nabla^W 
+\phi^* \sum_{k=1}^{\infty}\sum_{i+j=k-1}\{i;\nabla^W;j\}(\frac{F}{\nu})\\
&&\nonumber\qquad 
-\phi^* \circ 
ad(\frac{1}{\nu}(G_z-\phi^{-1*}(G_y)_z) )  \circ\exp[ad(\frac{F}{\nu})] \\
&=&\nonumber
\phi^*\circ\exp[ad(\frac{F}{\nu})]\circ\nabla^W 
+\phi^* \sum_{k=1}^{\infty}\sum_{i+j=k-1}\{i;\nabla^W;j\}(\frac{F}{\nu})\\
&&\nonumber\qquad
-\phi^*\circ\exp[ad(\frac{F}{\nu})](\frac{1}{\nu}(G_z-\phi^{-1*}(G_y)_z))\\
&&\nonumber\qquad\qquad
-\phi^*\Bigr(\sum_{k=1}^{\infty} \frac{1}{k!}\sum_{i+j=k-1} 
\{i;\frac{1}{\nu}(G_z-\phi^{-1*}(G_y)_z);j\}\frac{F}{\nu}\Bigr) \\
&=&
\nonumber
\phi^*\circ \exp[ad(\frac{F}{\nu})]\circ \nabla^W\\
&&\nonumber
\qquad + \phi^*\sum_{k=1}^{\infty} \sum_{i+j=k-1}\{i;\tilde{\nabla}^W;j\}
(\frac{F}{\nu})\\
&&\nonumber
\qquad\qquad-\phi^*\Bigr(\exp[ad(\frac{F}{\nu})]
\bigr(\frac{1}{\nu}(G_z-\phi^{-1*}(G_y)_z)\bigr)\Bigr). 
\end{eqnarray}
$\square$
\end{pf}
Therefore we have 
\begin{thm} 
Assume that $F$ satisfies 
\begin{equation}\label{lift-equation}
\sum_{k=1}^\infty \frac{1}{k!} \sum_{i+j=k-1} 
\{i;\tilde{ \nabla }^W;j\}(\frac{F}{\nu}\bigl)=
\exp[ad\bigl(\frac{F}{\nu}\bigl)]\bigl(\frac{1}{\nu}(G-\phi^{-1*}(G))\bigl)
\end{equation}
where $G$ and $\phi$ is given in (\ref{phi}) and  
(\ref{G}). 
Then we have  
\begin{equation}\nabla^W\circ\phi^*\circ \exp[ad(\frac{1}{\nu}F)]=
\phi^*\circ \exp[ad(\frac{1}{\nu}F)]\circ\nabla^W. 
\end{equation}
\end{thm}
With a slight modification, we can adapt the 
iterated argument employed 
for the construction of Fedosov connection (\cite{ f}, see also \S 6.1) 
in such a way that 
we can apply it to solving the equation (\ref{lift-equation}). 
Thus, we have 
\begin{thm}\label{lift-thm}
For any symplectic 
diffeomorphism $\phi$ on a symplectic manifold $(M,\omega)$, 
there canonically exists an element 
$\hat{\phi}\in Aut(M,*)$ which induces 
the base map $\phi$ on $M$. 
\end{thm}
\begin{prop}\label{proposition-bijection}
Assume that there exists a map 
\footnote{The map $j$ is not a Lie group homomorphism in general. }
$j$ from ${\rm Diff}(M,\omega)$ into ${\rm Aut}(M,*)$ satisfying 
$p\circ j={\rm identity}$. Then 
we have a bijection:
\begin{equation}\label{bijection}
 {\rm Aut}(M,*)\cong {\rm Diff}(M,\omega)\times \underline{\rm Aut}(M,*) . 
\end{equation}
\end{prop}
\begin{pf} 
As mentioned in Proposition~\ref{inducing}, 
any element $\Psi\in {\rm Aut}(M,*)$
induces a symplectic diffeomorphism $\phi=p(\Psi)$ on the base manifold. 
Set $\hat{\phi}=j(\phi)$ and $\Phi=\hat{\phi}^{-1}\circ\Psi$. 
By the assumption,   
$\Phi$ induces the base identity map. 
According to Propositions~\ref{representative2} and 
\ref{representative}, we see  
$\Phi=\exp[ad (\frac{1}{\nu}(g(\nu^2)+\nu^2f^\#(\nu^2)))]$. 
$\square$
\end{pf} 
As seen in the proposition above, in order to determine the 
model space of ${\rm Aut}(M,*)$, we have to determine the model space 
of ${\rm Diff}(M,\omega)$. Take a diffeomorphism $(pr_M,\sigma)$ 
from an open neighborhood $U_0$ of the zero section in $T^*M$ onto an 
open neighborhood $U_2$ of 
the diagonal set of $M\times M$, 
such that $\sigma(0\mbox{-section}|_x)=x$. 
Let $\omega_0$ be the canonical symplectic structure of 
$T^*M$, and $\omega_1:=(pr_M,\sigma)^*(\omega\oplus\omega^-)$, 
where the reversed symplectic structre of $\omega$ is denoted by 
$\omega^-$.
Since $\omega_0$ and $\omega_1$ vanish  
when restricted to the zero section, 
by vurtue of Moser's technique (cf. \cite{banyaga}), there exists a 
diffeomorphism $\varphi:U_0 \rightarrow U_1$ between 
two suitable 
open neighborhoods $U_0$ and $\,U_1$ of the zero section in $T^*M$ 
which is the identity on the zero section and 
satisfies $\varphi^*\omega_1=\omega_0$. Thus we obtain that  
\begin{equation}\label{diagram} 
\begin{array}{cccccc}
\eta=(pr_M,\sigma)\circ\varphi:& (U_0,\omega_0) & 
\stackrel{\varphi}{\longleftrightarrow} & (U_1,\omega_1) 
& \stackrel{(pr_M,\sigma)}{\longleftrightarrow} & (U_2,\omega\oplus\omega^-) . 
\end{array}
\end{equation}
\par\bigskip\noindent
\unitlength 0.1in
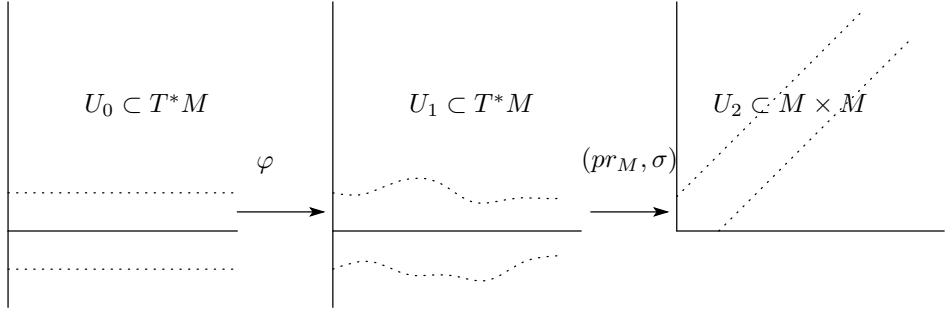
\begin{picture}( 51.0000, 16.0000)(  8.0000,-28.0000)
%
\special{pn 8}%
\special{pa 800 2400}%
\special{pa 2000 2400}%
\special{fp}%
%
\special{pn 8}%
\special{pa 2500 2400}%
\special{pa 3800 2400}%
\special{fp}%
%
\special{pn 8}%
\special{pa 2500 2400}%
\special{pa 2500 2800}%
\special{fp}%
%
\special{pn 8}%
\special{pa 800 2800}%
\special{pa 800 2400}%
\special{fp}%
%
\special{pn 8}%
\special{pa 800 2200}%
\special{pa 2000 2200}%
\special{dt 0.045}%
%
\special{pn 8}%
\special{pa 800 2600}%
\special{pa 2000 2600}%
\special{dt 0.045}%
%
\special{pn 8}%
\special{pa 2500 2200}%
\special{pa 2532 2206}%
\special{pa 2566 2210}%
\special{pa 2598 2212}%
\special{pa 2628 2212}%
\special{pa 2660 2208}%
\special{pa 2690 2200}%
\special{pa 2720 2188}%
\special{pa 2748 2176}%
\special{pa 2778 2164}%
\special{pa 2810 2152}%
\special{pa 2840 2140}%
\special{pa 2872 2132}%
\special{pa 2904 2126}%
\special{pa 2938 2124}%
\special{pa 2970 2124}%
\special{pa 3002 2128}%
\special{pa 3032 2138}%
\special{pa 3062 2150}%
\special{pa 3090 2166}%
\special{pa 3118 2184}%
\special{pa 3146 2204}%
\special{pa 3174 2222}%
\special{pa 3200 2236}%
\special{pa 3230 2248}%
\special{pa 3258 2254}%
\special{pa 3290 2254}%
\special{pa 3320 2252}%
\special{pa 3354 2246}%
\special{pa 3386 2240}%
\special{pa 3418 2234}%
\special{pa 3452 2230}%
\special{pa 3484 2226}%
\special{pa 3516 2226}%
\special{pa 3548 2228}%
\special{pa 3580 2228}%
\special{pa 3612 2230}%
\special{pa 3642 2230}%
\special{pa 3674 2230}%
\special{pa 3690 2230}%
\special{sp -0.045}%
%
\special{pn 8}%
\special{pa 2500 2600}%
\special{pa 2532 2588}%
\special{pa 2564 2578}%
\special{pa 2596 2570}%
\special{pa 2628 2562}%
\special{pa 2658 2558}%
\special{pa 2690 2560}%
\special{pa 2720 2564}%
\special{pa 2750 2574}%
\special{pa 2780 2586}%
\special{pa 2810 2598}%
\special{pa 2840 2612}%
\special{pa 2870 2622}%
\special{pa 2900 2630}%
\special{pa 2932 2634}%
\special{pa 2964 2634}%
\special{pa 2996 2632}%
\special{pa 3028 2630}%
\special{pa 3060 2630}%
\special{pa 3092 2630}%
\special{pa 3126 2634}%
\special{pa 3158 2640}%
\special{pa 3190 2648}%
\special{pa 3222 2656}%
\special{pa 3252 2660}%
\special{pa 3284 2662}%
\special{pa 3314 2658}%
\special{pa 3342 2648}%
\special{pa 3370 2636}%
\special{pa 3398 2618}%
\special{pa 3426 2602}%
\special{pa 3456 2584}%
\special{pa 3484 2568}%
\special{pa 3514 2556}%
\special{pa 3544 2546}%
\special{pa 3576 2540}%
\special{pa 3608 2536}%
\special{pa 3640 2532}%
\special{pa 3672 2530}%
\special{pa 3680 2530}%
\special{sp -0.045}%
%
\special{pn 8}%
\special{pa 2000 2300}%
\special{pa 2450 2300}%
\special{fp}%
\special{sh 1}%
\special{pa 2450 2300}%
\special{pa 2384 2280}%
\special{pa 2398 2300}%
\special{pa 2384 2320}%
\special{pa 2450 2300}%
\special{fp}%
%
\special{pn 8}%
\special{pa 3850 2300}%
\special{pa 4250 2300}%
\special{fp}%
\special{sh 1}%
\special{pa 4250 2300}%
\special{pa 4184 2280}%
\special{pa 4198 2300}%
\special{pa 4184 2320}%
\special{pa 4250 2300}%
\special{fp}%
\put(21.0000,-21.0000){\makebox(0,0)[lb]{$\varphi$}}%
\put(38.0000,-21.0000){\makebox(0,0)[lb]{$(pr_M,\sigma)$}}%
%
\special{pn 8}%
\special{pa 800 2400}%
\special{pa 800 1200}%
\special{fp}%
%
\special{pn 8}%
\special{pa 2500 2400}%
\special{pa 2500 1200}%
\special{fp}%
%
\special{pn 8}%
\special{pa 4300 2400}%
\special{pa 4300 1200}%
\special{fp}%
\put(12.0000,-18.0000){\makebox(0,0)[lb]{$U_0\subset T^*M$}}%
\put(29.0000,-18.0000){\makebox(0,0)[lb]{$U_1\subset T^*M$}}%
%
\special{pn 8}%
\special{pa 5900 2400}%
\special{pa 5900 2400}%
\special{fp}%
%
\special{pn 8}%
\special{pa 4300 2400}%
\special{pa 5700 2400}%
\special{fp}%
%
\special{pn 8}%
\special{pa 4300 2220}%
\special{pa 5270 1250}%
\special{dt 0.045}%
%
\special{pn 8}%
\special{pa 4520 2400}%
\special{pa 5520 1400}%
\special{dt 0.045}%
\put(44.9000,-18.0000){\makebox(0,0)[lb]{$U_2\subset M\times M$}}%
\end{picture}%
\par\medskip
\noindent
We also see that 
\begin{eqnarray}
\nonumber&&\{\eta^{-1}(x,f(x))|x\in M\}
\mbox{ is a closed form } (\in \Omega_c^1(T^*M) ) \\
\nonumber&\Leftrightarrow&\{\eta^{-1}(x,f(x))|x\in M\}
\mbox{ is a Lagrangian submanifold of }(T^*M,\omega_0) \\
&\stackrel{(\ref{diagram})}{\Leftrightarrow}&
\label{equiv}
\mbox{ the graph is a Lagrangian 
submanifold of}~(M\times M,\omega\oplus\omega^-) 
\nonumber
\\
\nonumber
&\Leftrightarrow&
0=(Id_M,f)^*(pr_1^*(\omega)-pr_2^*(\omega))=Id_M^*\omega-f^*\omega  
\\
\nonumber
&\Leftrightarrow& 
f\in {\rm Diff}_c(M,\omega)
\end{eqnarray}
Let $\frak{U}$ be an open neighborhood of $Id_M$ consisting of all 
$f\in {\rm Diff}(M)$ with compact support satisfying 
$(Id_M,f)(M)\subset U_2$ and $pr_M
:\eta^{-1}(\{(x,f(x))|x\in M\})\rightarrow M$ is 
still a diffeomorphism. 
For $f\in {\frak{U}}$, 
the map 
$(Id_M,f):M \rightarrow graph(f)\subset M\times M$ is 
the natural diffeomorphism onto the graph of $f$. 
According to (\ref{diagram}), we can 
define the smooth chart of ${\rm Diff}(M)$ 
which is centered at the identity in the following way: 
\begin{eqnarray}\label{chart}\nonumber
&&\nonumber
{\rm Diff}_c(M)\supset {\frak{U}} \stackrel{\Psi}{\rightarrow} \Psi({\frak{U}})
\subset 
\Omega_c^1(M), \quad 
\Psi(f)=\eta^{-1}(Id_M,f)~;~ M \rightarrow T^*M.
\end{eqnarray}
Since 
$\Omega_c^1(T^*M)$ is Mackey complete (cf. \cite{km}),  
${\frak{U}} \cap {\rm Diff}(M,\omega)$ 
gives a 
submanifold chart for ${\rm Diff}(M,\omega)$ at $Id_M$.   
Moreover, conditions of Definition~\ref{regular}
can be shown by the standard argument 
of ordinary differential equation 
under a certain identification 
of $T^\ast M$ with $TM$. 
Therefore, we have the following. 
\begin{thm}[\cite{km}, \cite{o}]\label{Lie-group} 
Let $(M,\omega)$ be a finite-dimensional symplectic manifold. 
Then the group ${\rm Diff}(M,\omega)$ of symplectic diffeomorphisms is 
a regular Lie 
group and a closed submanifold of the regular Lie group ${\rm Diff}(M)$ 
of diffeomorphisms. The Lie algebra of ${\rm Diff}(M,\omega)$ 
is Mackey complete locally convex space 
${\mathfrak X}_c(M,\omega)$ of 
symplectic vector fields with compact supports.  
\end{thm}
Combining the Baker-Campbell-Hausdorff formula with 
Propositions \ref{regularity-underline-aut}  
and \ref{proposition-bijection}, we have 
\begin{lem}
\label{smooth-operations}
The following maps are smooth: 
\begin{enumerate}
\item[]{\rm (i)} ${\rm Diff}(M,\omega)\times \underline{\rm Aut}(M,*) 
           \rightarrow 
                  \underline{\rm Aut}(M,*)
      ;(\phi,\Psi)
           \mapsto  
                   \hat{\phi}^{-1}\circ \Psi\circ\hat{\phi}$,
\item[]{\rm (ii)} ${\rm Diff}(M,\omega)\times {\rm Diff}(M,\omega)
            \rightarrow  
                 \underline{\rm Aut}(M,*)
      ;(\phi,\psi)
            \mapsto 
                 \widehat{(\phi\circ\psi)}^{-1}\circ\hat{\phi}\circ\hat{\psi}$, \item[]{\rm (iii)} ${\rm Diff}(M,\omega)\rightarrow \underline{\rm Aut}(M,*)
      ;\phi \mapsto \hat{\phi}\circ\hat{\phi^{-1}}$.
\end{enumerate}
\end{lem}
According to Propositions \ref{regularity-underline-aut}  
and \ref{proposition-bijection}, 
$\frak{X}_c(M,\omega)\times\frak{ C}_{c}(M)$ is a model space, 
which is a Mackey complete locally convex space.  
Let $\Psi_i=\hat{\psi_i}\circ e^{ad(\frac{1}{\nu}H_i(\nu^2))}$, where 
$H_i(\nu^2)=g_i(\nu^2)+\nu^2 f_i^\#(\nu^2)$ $(i=1,2).$
Then the multiplication is written in the following way: 
\begin{eqnarray}\label{compo}
\Psi_1\circ\Psi_2
&=&\hat{\psi_1}\circ e^{ad(\frac{1}{\nu}H_1(\nu^2))}\circ 
   \hat{\psi_2}\circ e^{ad(\frac{1}{\nu}H_2(\nu^2))}\\
\nonumber&=&\widehat{\psi_1\circ\psi_2}\circ\left\{
      \widehat{(\psi_1\circ\psi_2)}^{-1}
        \circ\hat{\psi_1}\circ\hat{\psi_2}\right\}\\
        &&\nonumber\qquad\circ
           \left\{ \hat{\psi_2}^{-1}
           \circ 
               e^{ad(\frac{1}{\nu}H_1(\nu^2))}\circ\hat{\psi_2}\right\}
                  \circ e^{ad(\frac{1}{\nu}H_2(\nu^2))}. 
\end{eqnarray}
According to {\rm (i)} and {\rm (ii)} of 
Lemma~\ref{smooth-operations}, $($\ref{compo}$)$ is written as 
$$ \left(\widehat{\psi_1\circ\psi_2}\right)\circ 
e^{ad(\frac{1}{\nu} H(\psi_1,\psi_2,H_1(\nu^2),H_2(\nu^2)))}, $$
and we see  
the smoothness of 
$$(\psi_1,\psi_2,H_1(\nu^2),H_2(\nu^2))\mapsto 
H(\psi_1,\psi_2,H_1(\nu^2),H_2(\nu^2)) . $$ 
By a similar way, we can verify the smoothness of the inverse operation. 
Summing up what is mentioned above, we have  
\begin{thm}\label{main2}
Under the assumption of Proposition~\ref{proposition-bijection}, 
$ {\rm Aut}(M,*)$ is a Lie group modeled on 
a Mackey complete locally convex space 
$\frak{X}_c(M,\omega)\times \frak{C}_c(M)$. 
\end{thm} 
Furthermore, combining the definition of ${\rm Aut}(M, *)$ with 
Proposition~\ref{inducing} gives a short exact sequence 
$$1\rightarrow \underline{\rm Aut}(M,*)\rightarrow {\rm Aut}(M,*) \rightarrow 
{\rm Diff}(M,\omega)\rightarrow 1.$$
As mentioned in Thenorem~\ref{Lie-group}, 
the group ${\rm Diff}(M,\omega)$ of all symplectic diffeomorphisms  
is a regular Lie group modeled on 
a Mackey complete locally convex space 
$\frak{X}_c(M,\omega)$. 
Therefore, combining Theorem~\ref{regularity-underline-aut} with 
Lemma~\ref{short-exact-sequence}, 
${\rm Aut}(M,*)$ is a regular Lie group. 
Thus, we obtain the following. 
\begin{thm}
\label{main3}
Under the same notation above, 
\begin{enumerate}
\item $1\rightarrow \underline{\rm Aut}(M,*) 
\rightarrow {\rm Aut}(M,*) \rightarrow 
{\rm Diff}(M,\omega)\rightarrow 1$
is a short exact sequence of Lie groups.   
\item 
$\underline{\rm Aut}(M,*)$ and ${\rm Aut}(M,*)$ are regular Lie groups. 
\end{enumerate}
\end{thm}
This completes the proof of Theorem~\ref{main-theorem-omori-proceedings}.

\subsection{Lifts as local contact Weyl diffeomorphisms}
\label{subsection-5-3}
We first remark that we can find a {\bf globally defined 
modified contact Weyl diffeomorphism} as 
a lift of symplectic diffeomorphism. 
However, in general, we can not find a {\bf globally defined 
contact Weyl diffeomrphism} as a lift. 
In the present subsection, we consider the existence of 
a {\bf locally defined} contact Weyl diffeomorpshism as a lift of 
a locally defined symplectic diffeomorphism. 
Although the following argument  
seems well known for specialists, we review it for readers' convenience. 

Assume that  
$$\bigr(U,z=(z^1,\cdots,z^{2n})\bigr), \quad
\bigr(\phi(U),z'=(z^{'1},\cdots,z^{'2n})\bigr)$$ 
are star-shaped Darboux charts.  
Then $\phi|_U$ is expressed as 
$$\bigr(z^{'1},\cdots,z^{'2n}\bigr) = 
\bigr(\phi^1(z),\cdots,\phi^{2n}(z)\bigr)$$ 
and satisfies 
$$\{\phi^i,\phi^j\}=\{\phi^{i+n},\phi^{j+n}\}=0,~~
\{\phi^i,\phi^{n+j}\}=-\delta^{ij}\qquad (1\le i,j\le n) , $$
because $\phi$ is a symplectic diffeomorphism defined on $U$. 
The Weyl continuations $\phi^{i\#}~(i=1,\cdots 2n)$ 
only satisfy
$$\!\!\!\!\!
[\phi^{i \#},\phi^{j \#}]
=
\nu^3 a_{(3)}^{i,j\#}+\cdots+\nu^{2l+1}a_{(2l+1)}^{i,j\#}+\cdots,
$$
\begin{equation}
\label{omy-global-3.1}
\qquad\qquad [\phi^{i \#},\phi^{n+j \#}]
=
-\nu\delta^{ij}
+\nu^3 a_{(3)}^{i,n+j \#}+\cdots+\nu^{2l+1} a_{(2l+1)}^{i,n+j \#}+\cdots, 
\end{equation}  
$$\qquad\qquad\!
[\phi^{n+i \#},\phi^{n+j \#}]
=
\nu^3 a_{(3)}^{n+i,n+j \#}+\cdots+\nu^{2l+1} a_{(2l+1)}^{n+i,n+j \#}+\cdots.  
$$
However the Jacobi identity holds: 
\begin{equation}\label{jacobi}
[\phi^{s\#}[\phi^{t\#},\phi^{u\#}]]+\mbox{c.p.}=0,
\end{equation}
where ``c.p." means ``cyclic permutation".
This gives  
\begin{equation}\label{omy-global-3.2}
\{z^{'s},a_{(3)}^{t,u}\}+\mbox{c.p.}= 
\{\phi^s,a_{(3)}^{t,u}\}+\mbox{c.p.}=0\qquad (1\le s,t,u\le 2n). 
\end{equation}
Set 
\begin{eqnarray}
\omega'(z')&\nonumber=& \frac{1}{2}\sum_{1\le i,j\le n}\Bigr[
a_{n+i,n+j}^{(3)} (z')dx^{'i}\land dx^{'j}\\
&&\qquad -2 a_{n+i,j}^{(3)}(z') dx^{'i}\land dy^{'j} +a_{i,j}^{(3)} (z') 
dy^{'i}\land dy^{'j}\Bigr]\quad(z'\in U'). 
\end{eqnarray}
A direct computation shows that 
(\ref{omy-global-3.2}) is equivalent to $d \omega'=0$. 
As the proof of Lemma~3.4 in \cite{omy2}, 
the closedness of $\omega'$ above ensures the 
existence of elements 
$b'_{j}
\in C^\infty(\phi(U))[[\nu]],~~(j=1,\cdots, 2n)$ 
such that  
replacing $\phi^{s\#}$ by
\begin{equation}
\label{omy-global-3.5}
\phi_{(1)}^s=
\left\{
\begin{array}{ll}
\phi^j(z)+\nu^2b'_{j+n}(\phi(z)), & s=j \\ 
\phi^{j+n}(z)-\nu^2b'_j(\phi(z)),  & s=j+n 
\end{array}\right.
\qquad (1\le j \le n),  
\end{equation}
shows that 
$\nu^3$-components of 
(\ref{omy-global-3.1}) vanish. 
Repeating the argument above for the $\nu^5$-, $\nu^7$- components 
gives 
$$\phi_{(\infty)}=(\phi_{(\infty)}^1,\cdots,\phi_{(\infty)}^{2n}) , $$ 
where 
\begin{equation}
\label{i}
\phi_{(\infty)}^i=\phi^i(z)+\sum_{p\ge 1} \nu^{2p} g_p^i(z)
\end{equation} 
such that
$$[\phi_{(\infty)}^{i\#},\phi_{(\infty)}^{j\#}]=
  [\phi_{(\infty)}^{n+i\#},\phi_{(\infty)}^{n+j\#}]=0, 
\quad  [\phi_{(\infty)}^{i\#},\phi_{(\infty)}^{n+j\#}]=-\nu\delta^{ij},~
(i,j=1,\cdots, n) . 
$$
Thus, by Lemma 3.2 in \cite{omy}, 
there exists a local Weyl diffeomorphism 
$\Phi_U$ which induces the base map $\phi_U$. 
We next extend $\Phi_U$ to a local contact Weyl diffeomorphism
$\Psi_U$. 
Set 
\begin{equation}\label{ii}
\Psi_U^*(a)=\left\{\begin{array}{ll}
\Phi^*_U(a),&(a\in {\mathcal F}_U), \\
\tilde{\tau}_U+H,&(a=\tilde{\tau}_{\phi(U)}).
\end{array}
\right.
\end{equation}
where $H=\sum_m\nu^mh_m^\#$ is an unknown term. 
$\Psi_U$ is a contact Weyl diffeomorphism if it satisfies the following 
equation w.r.t. $H$  
\begin{equation}\label{iii}
[\Psi_U^* (\tilde{\tau}_{\phi(U)} ), \Psi_U^*(z^{'i\#}) ]=\Psi_U^*
[\tilde{\tau}_{\phi(U)}, z^{'i\#}]
. 
\end{equation}
As to the equation, 
we easily have 
\begin{equation}
\mbox{R.H.S. of } ( \ref{iii})=
\Psi_U^*(\nu z^{'i\#})
\stackrel{def}{=}\nu ( \phi_i^{\#}+B^{\#}(\nu)),
\label{iv}
\end{equation}
where $B(\nu)=\sum_{l\ge 1} \nu^{2l}g_{l}$. 
On the other hand, we also obtain 
\begin{eqnarray}
\label{v}
&\nonumber&\mbox{L.H.S. of }(\ref{iii})\\
&{\stackrel{(2.18)\mbox{ in }\cite{y}}{=}}&
{\nu}{\Bigr(}\sum_lz^l\frac{\partial z^{'i}}{\partial z^l}{\Bigr)}^\#
{+}{\Bigr[}\sum_m\nu^mh_m(z^{'i}\circ\phi){+}\sum_p\nu^{2p}g_p{\Bigr]}^{\#}
\\ 
&\nonumber& \qquad\qquad 
{+}{\Bigr(}2\nu^2\partial_\nu B {+}\nu (EB){\Bigr)}^{\#}
\end{eqnarray}
where $E= \nu \sum_{l=1}^{2n}z^l\partial_{z^l}$. 
As the proof of Theorem~3.6 in \cite{omy}, comparing the components 
w.r.t. $\nu^1$-, $\nu^2$-,$\nu^3$-,$\cdots $ 
of the both sides splits the equation w.r.t. $H$ above 
into infinitely many equations.  
Since the component of $\nu$ is 
\begin{equation}
\{h_0,z^{'i}\circ\phi\}=(z^{'i}\circ\phi)
-\sum z^l\left(\frac{\partial z^{'i}}{\partial z^l}\right)
,
\end{equation}
we can find the solution $h_0$ for this equation, 
and then we can solve the infinitely many equations recursively
\footnote{Thanks to star-shapedness of $U$, we 
can fix $b'_{s}$ and $H$ canonically. }. 
Summing up the above, we have 
\begin{prop}\label{contactweyllift}
Take a star-shaped Darboux chart $U$.  
For any symplectic diffeomorphism 
$\phi:U\rightarrow \phi(U)$, 
there canonically exists a contact Weyl  
diffeomorphism~{\rm (}CW-lift{\rm )} $\hat{\phi}$ which induces $\phi$. 
\end{prop}
Then we have 
\begin{cor}
Assume that a symplectic manifold $M$ is covered by a star-shaped 
Darboux chart. 
Then for any symplectic diffeomorphism 
$\phi:M \rightarrow M$, 
there canonically exists a contact Weyl  
diffeomorphism~{\rm (}CW-lift{\rm )} $\hat{\phi}$ which induces $\phi$. 
\end{cor}

\section{Appendices}
\subsection{Fedosov connection}
As seen in the previous section, 
as to the quantum connection $\nabla^Q$, 
it holds that 
$$\begin{array}{l}
\nabla^Q|_{{W}_M} = \nabla^W, \\   
{\cal F}({W}_M)=\{\mbox{ parallel section w.r.t.}
 \nabla^Q|_{{W}_M}\}. 
\end{array} 
$$ 
Let $\nabla^{symp}$ be a symplectic connection and 
\begin{equation}
\delta^{-1}(\nu^mZ^\alpha dz^\beta)=
\left\{\begin{array}{l}
\sum_{i=1}^{2n} dz_i 
\iota_{Z_i}\nu^mZ^\alpha dz^\beta\quad(|\alpha|+|\beta|\not=0) ,\\
0  \quad (|\alpha|+|\beta|=0) , 
\end{array}\right. 
\end{equation}
where $\iota$ is a inner product. 
We may write $\nabla^Q|_{{W}_M}=\nabla^{symp}-\delta+r$, 
where $r$ is a 1-form with $\Gamma(W_M)$ coefficient. 
Then as in \cite{f}, $r$ satisfies the following equation 
\begin{equation}
\delta r = R_\omega + \nabla^{symp} r +\frac{1}{2\nu}[r,r]. 
\end{equation}
Or equivalently 
$r$ satisfies 
\begin{equation}\label{F-eq}
r=\delta^{-1}\{(\nabla^{symp} + \frac{1}{2\nu}[r,r])+R_\omega\}, 
\end{equation} 
under the assumptions 
$\mbox{deg}~r\ge 2, ~ 
\delta^{-1}r=0, ~ 
r_0=0$. 
Set $r_k$ is the term of $r $ degree $k$. 
Since it is known that this equation can be solved by recursively 
in the following way 
\begin{equation}\begin{array}{l}
r_{3}= \delta^{-1} R_\omega,\\ 
 r_{n+3}=\delta^{-1}(\nabla^Kr_{n+2}+\frac{1}{\nu}\sum_{l=1}^{k-1} 
r_{l+2}\ast r_{k+2-l}). 
\end{array}
\end{equation}
Until now we did not consider the symplectic action of $G$ on $M$. 
By the same manner, we have the following. 
\begin{prop}\label{commyoshiokaconn}
Suppose that a Lie group $G$ is compact, and 
for any $g\in G$, $\nabla^{symp} \circ g^\ast= 
g^\ast \circ \nabla^{symp}$, where $\nabla^{symp}$ is a symplectic 
connection. 
Then we can construct a quantum connection 
$\nabla^Q$ such that 
$\nabla^Q|_{{W}_M}\circ g^\ast=
g^\ast\circ \nabla^Q|_{{W}_M}$\footnote{Note that 
this action is not commute with $\nabla^Q$ in general. 
For example, compute and compare 
$\nabla^Q (g^\ast \tau)$ and $g^\ast(\nabla^Q\tau)$}. 
\end{prop}

\subsection{Examples of star exponential}
This subsection is devoted to computations of star-exponential functions for 
quadratic forms~(cf. \cite{ommy3}, \cite{ommy4}). 

Let $Z={}^t(Z^1,\ldots,Z^{2n})$, $A[Z]:={}^tZAZ$, where 
$A \in Sym(2n,{\mathbf R})$, i.e. 
$A$ is a $2n\times 2n$-real symmetric matrix. In order to compute 
the star exponential function with respect to the 
Moyal product
$e_*^{\frac{1}{\mu}A[Z]}$, 
we treat the following evolution equation. 
\begin{equation}\label{evolution}
\partial_t F=\frac{1}{\mu}A[Z]*F,\quad F_0=e^{\frac{1}{\mu}B[Z]}, 
\end{equation}
where $B\in Sym(2n,{\mathbf R}),~\mu=-\sqrt{-1}\hbar.$ 
Under the assumption $F(t)=g\cdot e^{\frac{1}{\mu}Q[Z]}$ 
($g=g(t),~Q=Q(t)$), we would like to find a solution of this equation. 
Set $\Lambda=\left[\begin{array}{cc}0 & 1 \\-1 & 0 \end{array}\right]$, 
$q:=\Lambda Q$ and $a:=\Lambda A$, 
then we see 
that 
\begin{eqnarray}\label{order}
&&
\sum_{l,m,i_1,i_2,j_1,j_2}
A_{i_1i_2}\Lambda^{i_1j_1}\Lambda^{i_2j_2}
Q_{j_1m}Q_{j_2l}Z^mZ^l\\
&=&
\nonumber
\sum_{l,m,i_1,i_2,j_1,j_2}
Q_{mj_1}(-\Lambda^{j_1i_1})A{i_1i_2}\Lambda^{i_2j_2}Q_{j_2l}Z^mZ^l\\
\nonumber
&=&-Q\Lambda A \Lambda Q[Z]. 
\end{eqnarray}
From now, we use Einstein's convention. 
As to the first equation of (\ref{evolution}), 
we see that 
\begin{eqnarray}\label{LHS}
\nonumber
\mbox{L.H.S. of } (\ref{evolution}) &=& 
g' e^{\frac{1}{\mu}Q[Z]} +g {\frac{1}{\mu}Q'[Z]}e^{\frac{1}{\mu}Q[Z]}, \\
\label{RHS}
\nonumber
\mbox{R.H.S. of } (\ref{evolution}) &=& {\frac{1}{\mu}A[Z]}*F  \\
\nonumber&
{=}&
{\frac{1}{\mu}A[Z]}\cdot F+\frac{i\hbar}{2}
\Lambda^{i_1j_1}\partial_{i_1}{\frac{1}{\mu}A[Z]}\cdot \partial_{j_1}F
-\frac{\hbar^2}{2\cdot 4} \Lambda^{i_1j_1}
\Lambda^{i_2j_2} \partial_{i_1i_2}{\frac{1}{\mu}A[Z]}
\partial_{j_1j_2}F\\
\nonumber
&=&
\frac{1}{\mu}A[Z]\cdot g e^{\frac{1}{\mu}Q[Z]}
- \frac{\mu}{2} \Lambda^{i_1j_1}
\Bigr(\frac{2}{\mu}A_{i_1l}Z^l\Bigr)
\Bigr(2g \frac{1}{\mu} Q_{j_1m}Z^me^{\frac{1}{\mu}Q[Z]}\Bigr)+ \\ 
\nonumber&& 
\qquad + \frac{\mu^2}{8}\Lambda^{i_1j_1}
\Lambda^{i_2j_2}
\Bigr(\frac{2}{\mu}A_{i_1i_2}\Bigr)\times\\
&&\nonumber\qquad\qquad\times 
\Bigr(2g \frac{1}{\mu} Q_{j_1j_2} e^{\frac{1}{\mu}Q[Z]}
+4g \frac{1}{\mu} Q_{j_1m}\frac{1}{\mu}*Q_{j_2l} Z^mZ^l e^{\frac{1}{\mu}Q[Z]} 
\Bigr) . 
\end{eqnarray}
Comparing the coefficient of $\mu^{-1}$, we obtain 
\begin{equation}
Q'[Z]=A[Z]-2{}^tA\Lambda Q[Z]-Q\Lambda A \Lambda Q[Z] . 
\end{equation}
Applying $\Lambda$ by left, we get 
\begin{eqnarray}
\Lambda Q' &=& \Lambda A + 
\Lambda Q\Lambda A-\Lambda A\Lambda Q 
-\Lambda Q\Lambda A\Lambda Q \\
\nonumber&=&
(1+\Lambda Q)\Lambda A(1-\Lambda Q)\\
\nonumber&=&(1+q)a(1-q).
\end{eqnarray}
As to the coefficient of $\mu^0$, we have 
\begin{eqnarray}
g'&=&\frac{1}{2}\Lambda^{i_1j_1}\Lambda^{i_2j_2}A_{i_1i_2}gQ_{j_1j_2} \\ 
\nonumber&=&-\frac{1}{2}tr(\Lambda A\cdot\Lambda Q)g\\
\nonumber&=&-\frac{1}{2}tr(aq) \cdot g . 
\end{eqnarray}
Thus the equation (\ref{evolution}) is rewritten by 
\begin{eqnarray}
\partial_t q &=& (1+q)a(1-q), \\ 
\partial_t g &=& -\frac{1}{2}tr(aq)\cdot g. 
\end{eqnarray}
We now recall the ``Cayley transform." 
\begin{prop}\label{cayley}
Set $C(X):=\frac{1-X}{1+X}$. Then 
\begin{enumerate}
\item $X\in sp(n,{\mathbf R})\Longleftrightarrow
\Lambda X\in Sym(2n,{\mathbf R})$, \\ and then 
$C(X)\in Sp(n,{\mathbf R}):=
\{g\in M(2n,{\mathbf R})|{}^tg\Lambda g=\Lambda\}$,  
\item $C^{-1}(g)=\frac{1-g}{1+g}$, 
\item $e^{2\sqrt{-1}a}=c(-\sqrt{-1}\tan(a))$, 
\item $\log a = 2\sqrt{-1}\arctan (\sqrt{-1}C^{-1}(g))$, 
\item $\partial_t q =(1+q)a(1-q)\Longleftrightarrow 
\partial_t C(q)=-2aC(q).$ 
\end{enumerate}
\end{prop}
\begin{pf} 
First we remark that 
\begin{eqnarray}
&&
(1-{}^tX)\Lambda(1-X)-(1+{}^tX)\Lambda(1+X)\\
\nonumber&=&\Lambda -{}^tX\Lambda-\Lambda X+{}^tX\Lambda X-
(\Lambda -{}^tX\Lambda +\Lambda X+{}^tX\Lambda X)=0, 
\end{eqnarray}
if ${}^tX\Lambda +\Lambda X=0$. 
Hence 
\begin{eqnarray}\label{star-key}
&&{}^tC(X)\Lambda C(X) =\Lambda\\
\nonumber
&=&{}^t\Bigr(\frac{1-X}{1+X}\Bigr)\Lambda\Bigr(\frac{1-X}{1+X}\Bigr)\\
\nonumber
&=& (1+{}^tX)^{-1}(1-{}^tX)\Lambda(1-X)(1+X)^{-1}. 
\end{eqnarray}
Take an element $A\in Sym(2n,{\mathbf R})$, and set
$X=\Lambda A$. 
Then 
$${}^t(\Lambda A)\Lambda+\Lambda(\Lambda A)\stackrel{(\ref{star-key})}{=}0.$$ 
Thus, $X=\Lambda A\in Sp(n,{\mathbf R})$. 
Conversely, assume 
$X\in Sp(n,{\mathbf R})$. Then 
$${}^t(-\Lambda X)=-{}^tX{}^t\Lambda ={}^t\Lambda=-\Lambda X. $$ 

2 is obvious. 

As to the assertion 3, 
\begin{eqnarray}
C^{-1}(e^{2\sqrt{-1}a})&=&\frac{1-e^{2\sqrt{-1}a}}{1+e^{2\sqrt{-1}a}}\\
&=&\nonumber
\frac{-\sqrt{-1}\frac{e^{\sqrt{-1}a}-e^{-\sqrt{-1}a}}{2\sqrt{-1}}}
{\frac{e^{\sqrt{-1}a}+e^{-\sqrt{-1}a}}{2}}\\
&=&\nonumber -\sqrt{-1}\tan a. 
\end{eqnarray}

As to the assertion 4, 
according to the assersion 3, for $g=e^{2\sqrt{-1}a}$, we easily have 
$$g=C(-\sqrt{-1}\tan \bigr(\frac{1}{2\sqrt{-1}}\log g\bigr)).$$
Then we see that 
$$\log g = 2\sqrt{-1} \arctan  (\sqrt{-1}C^{-1}(g)). $$

Finally we show the assertion 5. 
\begin{eqnarray}\label{evolution-Cayley-0}
C(q)'&=&  
\Bigr(\frac{1-q}{1+q}\Bigr)'\\
\nonumber&=&(1+q)^{-1}(-q)'+(1+q)^{-1}(-q)'(1+q)^{-1}(1-q)\\
\nonumber&=&-(1+q)^{-1}\{(1+q)a(1-q)\}
\\
\nonumber&&\qquad +(1+q)^{-1}\{(1+q)a(1-q)\}(1+q)^{-1}(1-q)\\
\nonumber&=&
-a(1-q)-a(1-q)(1+q)^{-1}(1-q)\\
\nonumber&=&-a\Bigr\{1+\frac{1-q}{1+q}\Bigr\}(1-q)\\
&=&\nonumber-2aC(q).
\end{eqnarray}
\end{pf}
Solving the above equation (\ref{evolution-Cayley-0}), 
we have 
$$C(q)=e^{-2at}C(b),$$ 
and then 
$$
q=C^{-1}(e^{-2at}\cdot C(b))=C^{-1}(C(-\sqrt{-1}\tan(\sqrt{-1}at)\cdot C(b )).
$$
Thus we obtain 
\begin{equation}\label{Q}
Q=-\Lambda\cdot C^{-1}(C(-\sqrt{-1}\tan (\sqrt{-1} \Lambda A t))\cdot C(b)). 
\end{equation} 
We can get $q$ in the following way. 
\begin{eqnarray}
\nonumber
q&=&(1-e^{-2at}C(b))(1+e^{-2at}C(b))^{-1}\\
\nonumber
&=&\Bigr(1-e^{-2at}\frac{1-b}{1+b}\Bigr)\Bigr(1-e^{-2at}\frac{1-b}{1+b}\Bigr)\\
&=&\nonumber
e^{-at}
\{ \}
(1+b)^{-1}\{(1+b)^{-1}\}^{-1} 
\{ \}
e^{-at}\\
&=&\nonumber
e^{-at}
\{e^{at}(1+b) -e^{-at}(1-b)\}
(1+b)^{-1}\\
&&\nonumber\qquad\times\{(1+b)^{-1}\}^{-1} 
\{e^{at}(1+b)+e^{-at}(1-b) \}^{-1}
(e^{-at})^{-1}
\\
&=&\nonumber
e^{-at}
\{e^{at}(1+b) -e^{-at}(1-b)\}
\{e^{at}(1+b)+e^{-at}(1-b) \}^{-1}
(e^{-at})^{-1} .  
\end{eqnarray}
This determines the phase part $Q$. 
Next we compute the amplitude coefficient part $g$. 

First we replace 
$$g'=-\frac{1}{2}Tr(aq)\cdot g$$
by 
\begin{equation}\label{log-diff-eq}
(\log g)'=-\frac{1}{2}Tr(aq).
\end{equation}
Since 
\begin{eqnarray}
&&\nonumber
Tr\Bigr\{\log \Bigr(\frac{e^{at}(1+b)+e^{-at}(1-b)}{2}
\Bigr)\Bigr\} '  
\\
&=&\nonumber
Tr\Bigr\{a\frac{e^{at}(1+b)-e^{-at}(1-b)}{e^{at}(1+b)+e^{-at}(1-b)}\Bigr\}
\\
&=&\nonumber
Tr(aq) , 
\end{eqnarray}
we can rewrite (\ref{log-diff-eq}) as 
\begin{eqnarray}
(\log g)'&=& 
-\frac{1}{2}
Tr\Bigr\{\log \Bigr(\frac{e^{at}(1+b)+e^{-at}(1-b)}{2}\Bigr)\Bigr\}'
\\
\nonumber
&=&-\frac{1}{2}\log\Bigr\{
\det \Bigr(\frac{e^{at}(1+b)+e^{-at}(1-b)}{2}\Bigr)\Bigr\}'.
\end{eqnarray}
Then we have 
$$g=\det{}^{-\frac{1}{2}}
\Bigr(\frac{e^{at}(1+b)+e^{-at}(1-b)}{2}\Bigr). $$ 
Setting $t=1$, $a=\Lambda A$ and $b=\Lambda B$, we get 
\begin{eqnarray}
&&\nonumber
e_*^{\frac{1}{\mu}A[Z]}*e^{\frac{1}{\mu}B[Z]}
\\
&=&\nonumber 
\det{}^{-\frac{1}{2}}\Bigr(\frac{e^{\Lambda A}(1+\Lambda B)
+e^{-\Lambda A}(1-\Lambda B) }{2}\Bigr)
\cdot 
e^{\frac{1}{\mu}\Lambda^{-1}C^{-1}[C(\frac{1}{\sqrt{-1}} 
\tan(\sqrt{-1}\Lambda A))\cdot C(\Lambda B)][Z]   } . 
\end{eqnarray}
Setting $B=0$, we have 
\begin{thm}
\begin{equation}\label{star_exponential}
e_*^{\frac{1}{\mu}A[Z]}
=\det{}^{-\frac{1}{2}}\Bigr(\frac{e^{\Lambda A}+e^{-\Lambda A}}{2}\Bigr)
\cdot e^{\frac{1}{\mu}
(\frac{\Lambda^{-1}}{\sqrt{-1}} \tan (\sqrt{-1}\Lambda A )) [Z]  } .
\end{equation}
\end{thm}

\subsection{Formality theorem} 
In this subsection, we recall the basics of $L_\infty$-algebras.
See \cite{k}, \cite{ds} and \cite{mk} for details.
\par
In the following $V=\oplus_{k\in {\mathbb Z}}V^k $ is a
graded vector space,
and $[1]$ is the shift-functor, that is, $V[1]^k=V^{k+1}$.
$V[1]=\oplus_k V[1]^{k}$ is called a shifted graded vector space of $V$.
We set
$C(V)=\oplus_{n\ge 1}Sym^n(V)$
where
$$
Sym^n(V)=T^n(V)/\{\cdots\otimes(x_1x_2-(-1)^{k_1k_2}x_2x_1)\otimes\cdots;
x_i\in V^{k_i}\} .
$$
This space has a coproduct
$\Delta:C(V)\rightarrow C(V)\otimes C(V)$ defined
in the following way:
\begin{eqnarray}\label{3}
&&\nonumber
\Delta(x_1\cdots x_n) \\
&\nonumber=&\sum_{k=1}^{n-1}\frac{1}{k!(n-k)!}
\sum_{\sigma\in S_n}{\rm sign}
(\sigma;x_1\cdots x_n)\bigr(x_{\sigma(1)}\cdots x_{\sigma(k)}\bigr) 
\otimes \bigr(x_{{\sigma}(k+1)}\cdots x_{\sigma( n)}\bigr) ,
\end{eqnarray}
where ${\rm sign}(\sigma;x_1\cdots x_n)$ is defined by
$
x_{\sigma(1)}\cdots x_{\sigma(n)}=
{\rm sign}(\sigma;x_1\cdots x_n)x_1\cdots  x_n.
$
This coproduct is coassociative, i.e.
$(1\otimes\Delta)\circ\Delta=(\Delta \otimes 1)\circ
\Delta .$
We denote $ k_1+k_2+\cdots +k_n$ by $\deg(x_1\cdots x_n) $,
 where ($x_i\in V^{k_i}$).
\begin{defn}
A map $f:C(V_1)\rightarrow C(V_2)$ is called a {\rm 
coalgebra homomorphism} if
{\rm (1)} $\Delta\circ f=(f\otimes f)\circ \Delta$,
{\rm (2)} $f$ preserves the grading.
\end{defn}
\noindent The coderivation is defined in the following way.
\begin{defn}
A map $\ell:C(V)\rightarrow C(V)$ is called a {\rm coderivation} if the
following properties are satisfied
{\rm (1)} $\ell$ is an odd vector field of degree $+1$,
{\rm (2)} $(\ell\hat{\otimes}
id+id\hat{\otimes}\ell)\circ \Delta=\Delta\circ \ell,$ where
$
(id\hat{\otimes}\ell)(x\otimes y)=(-1)^{\deg x}x\otimes \ell(y).
$
\end{defn}
\noindent We also use the following notation:
Set $f^{(n)}=p\circ f|_{Sym^n(V_1)}:{Sym^n(V_1)}\rightarrow V_2$, and
$\ell^{(n)}=p\circ \ell|_{Sym^n(V_1)}:{Sym^n(V_1)}\rightarrow V_2$,
where $p={\rm canonical~ projection}:C(V_2)\rightarrow V_2$.

Under the above notation, $L_\infty$-algebras and
$L_\infty$-morphisms are defined in the following way:
\begin{defn}
An $L_\infty$-{\rm algebra} is a pair $(V,\ell)$, where $V$ is a graded vector
space
and $\ell$ is a coderivation on the graded coalgebra $C(V)$, such that
$\ell^2=0$.
\end{defn}
\begin{defn}
An $L_\infty$-{\rm morphism} 
$F_*$ between two $L_\infty$-algebras $(V_1,\ell_1)$
and $(V_2,\ell_2)$ is a coalgebra homomorphism such
that $\ell_2\circ F_*=F_*\circ\ell_1$. \label{Linfty}
\end{defn}
\noindent{\bf Remark} If $\ell=\ell^{(1)}+\ell^{(2)}$, and
$d=\ell^{(1)},\,\,
[x,y]=(-1)^{\deg x -1} \ell^{(2)}(x,y),$ then
$\ell^2=0$ if and only if
\begin{eqnarray}
&&d^2=0\nonumber,
~~~d[x,y]=[dx,y]+(-1)^{\deg x-1}[x,dy], \\
&\nonumber &[[x,y],z]+(-1)^{(x+y)(z+1)}[[z,x],y]
+(-1)^{(y+z)(x+1)}[[y,z],x]=0, \nonumber
\end{eqnarray}
that is, $(V,\ell)$ is a graded differential Lie algebra.
\par\medskip
We next recall the Kontsevich formality theorem \cite{k}.
\par\bigskip
\noindent{\bf  Differential Graded Lie algebra of $T_{poly}$-fields}\par
\noindent 
Let $M$ be a smooth manifold. 
Set $T_{poly}(M)=\oplus_{k\ge -1}\Gamma(M,\land^{k+1}TM)$, and let 
$[\cdot,\cdot]_S$ be the {\it Schouten} bracket:
$$[X_0\land\cdots\land X_m,Y_0\land\cdots\land Y_n]_S=
\sum_{i,j}(-1)^{i+j+m}[X_i,Y_j]\cdots\land\hat{X_i}\land\cdots
\land\hat{Y_j}\land\cdots,
$$
$\mbox{where }X_i,Y_i\in \Gamma(M,TM).$ 
Then,  the triple 
$$(T_{poly}(M)[[\hbar]],d:=0,[\cdot,\cdot]:=[\cdot,\cdot]_S)$$ 
forms a differential graded Lie 
algebra. 
It is well known that for 
any bivector $\pi \in \Gamma(M,\land^2TM)$, 
$\pi$ 
is a Poisson structure if and only if 
\begin{equation}\label{sharp}
[\pi,\pi]_S=0.
\end{equation}
\par\medskip
\noindent{\bf  Differential Grade Lie algebra of $D_{poly}$-fields}\par
\noindent 
Let $(A,\bullet)$ be an associative algebra 
and set 
$
C(A)=\oplus_{k\ge -1}C^k$, $C^k=Hom(A^{\otimes k+1};A).
$ 
For $\varphi_i\in C^{k_i}$ $(i=1,2)$, we set 
\begin{eqnarray}
\label{14}
&&\nonumber\varphi_1\circ\hat{\varphi}_2(a_0\otimes a_1\otimes \cdots
\otimes a_{k_1+k_2})\\
&=& \sum_{i=0}^k(-1)^{ik_2}\varphi_1\left(a_0\otimes \cdots \otimes a_{i-1}
\right.
\\&&\qquad~~\nonumber\otimes\varphi_2 
\left.
(a_i\otimes \cdots \otimes a_{i+k_2})\otimes a_{i+k_2+1}\otimes \cdots \otimes a_{k_1+k_2}\right). 
\end{eqnarray}
Then the {\it Gerstenhaber} bracket is defined in the following way: 
\begin{equation}
\label{15}
[\varphi_1,\varphi_2]_G=\varphi_1\circ\hat{\varphi}_2-
(-1)^{k_1k_2}\varphi_2\circ\hat{\varphi}_1  
\end{equation}
and Hochschild coboundary operator $\delta=\delta_\bullet$ 
with respect to $\bullet$ is defined by 
$ 
\delta_\bullet(\varphi)=(-1)^k[\bullet,\varphi]
~~(\varphi\in C^k).
$ 
Then it is known that the triple 
$$(C(A),d:=\delta_\bullet,[\cdot,\cdot]:=[\cdot,\cdot]_G)$$
is a differential graded Lie algebra.
\par\medskip
Let $M$ be a smooth manifold. 
Set ${\cal F}=C^\infty(M)$, and 
$D_{poly}(M)^{n}(M)$ equals a space of all multidifferential 
operators from ${\cal F}^{\otimes n+1}$ into ${\cal F}$. 
Then $D_{poly}(M)[[\hbar]]=
\oplus_{n\ge -1}D_{poly}^{n}(M)[[\hbar]]$ is a subcomplex of $C({\cal F}
[[\hbar]])$.  
Furthermore,  the triple 
$(D_{poly}(M)[[\hbar]],\delta,[\cdot,\cdot]_G)$ 
is a differential graded Lie algebra. 
\par\medskip
\begin{prop}\label{hoch-ger}
Let $B$ be a bilinear operator and  
$f\star g = f\cdot g+B(f,g)$. 
Then the 
product 
$\star$ is associative if and only if $B$ satisfies 
\begin{equation}
\delta_\cdot B+\frac{1}{2}[B,B]_G=0.
\label{(7)} 
\end{equation}
\end{prop}
\par\bigskip\noindent
Next we recall the moduli space ${\cal MC}(C(V[1]))$. 
For $b\in V[1]$, set 
$e^b=1+b+\frac{b\otimes b}{2!}+\cdots\in C(V[1]). 
$ 
\begin{defn}
$
\ell(e^b)=0$  
is called a {\rm Batalin-Vilkovisky-Maurer-Cartan } equation, where 
$\ell =d+(-1)^{\mbox{deg }\circ}[\circ,\bullet]$. 
\label{BVMC}
\end{defn} Using this equation,  
we define the moduli space as follows:  
\begin{defn}
\begin{eqnarray}
\widehat{\cal MC}(C(V[1]))&=&\{b;\ell(e^b)=0 \} , \\
{\cal MC}(C(V[1]))&=&\widehat{\cal MC}(C(V[1]))/\thicksim, 
\end{eqnarray}
where $V$ stands for $T_{poly}(M)[[\hbar]]$ and $D_{poly}(M)[[\hbar]]$, and 
$\thicksim$ means the gauge equivalence
\footnote{Strictly speaking, as for formal Poisson bivectors, 
$\pi_1(\hbar)\sim\pi_2(\hbar)$ if 
there exists a formal vector field $D\in \frak{X}(M)[[\hbar]]$ such that 
$\exp \hbar D\circ \pi_1(\hbar)=
\pi_2(\hbar)\circ(\exp \hbar D\otimes\exp \hbar D)
$. 
On the other hand, as for star-prodcts, $*_1\sim *_2$ if 
there exists a intertwiner $T=1+\sum_{r\ge 1}\hbar^r
T_r,~(T_r\mbox{ :differential operators of order }r)$ such that $T\circ *_1=*_2\circ(T\otimes T ) $. 
} (cf. \cite{k}). 
\end{defn}
Note that  (\ref{sharp}) and (\ref{(7)}) can be seen as the  
Batalin-Vilkovisky-Maurer-Cartan equations. 
\par\medskip    
With these preliminaries, we can state precise version 
of Kontsevich formality theorem: 
\begin{thm} \label{Formality}
There exists a map ${\cal U}$ such that 
$${\cal U}:{\cal MC}(C(T_{poly}(M)[[\hbar]][1]))\cong{\cal MC}(C(D_{poly}(M)
[[\hbar]][1])).$$
\end{thm}
As a biproduct, we have 
\begin{thm}
For any Poisson manifold $(M,\omega)$ 
there exists a formal deformation quantization. 
\end{thm}

\par\vspace{.5cm}
\noindent
{\bf Acknowledgements}
The author is grateful 
to Professors Hideki Omori, Yoshiaki Maeda, Giuseppe Dito and 
Daniel Sternheimer 
for their advice and encouragement. 
He expresses his sincere gratitude to Professor Akira Yoshioka 
for fruitful discussions.

\end{document}